\documentclass[article]{elsarticle}
\usepackage[latin1]{inputenc}
\usepackage{amsfonts}
\usepackage{graphicx}
\usepackage{epsfig}
\usepackage{amssymb,amsmath,color,theorem}
\usepackage{amsmath}
\usepackage{subfigure}%

\newtheorem{thh}{Theorem}

\newtheorem{lemm}[thh]{Lemma}
\newtheorem{prop}[thh]{Proposition}

\newtheorem{definition}[thh]{Definition}
\newtheorem{remarque}{Remark}
\newtheorem{corollaire}[thh]{Corollary}

\journal{Journal of Computational and Applied Mathematics}

\graphicspath{
{E:/Multimedia/Images/Chercheurs/Dayan_Liu/JCAM/}
}

\begin{document}

\begin{frontmatter}

\title{Differentiation  by integration with Jacobi polynomials}

%\author{Da-yan Liu, Olivier Gibaru, Wilfrid
%Perruquetti}

\author[authorlabel1,authorlabel2]{Da-yan Liu}
\ead{dayan.liu@inria.fr}
\author[authorlabel1,authorlabel3]{Olivier Gibaru}
\ead{olivier.gibaru@ensam.eu}
\author[authorlabel1,authorlabel4]{Wilfrid
Perruquetti} \ead{wilfrid.perruquetti@inria.fr}

\address[authorlabel1]{\'{E}quipe Projet ALIEN, INRIA Lille-Nord Europe, Parc Scientifique
de la Haute Borne 40, avenue Halley B\^{a}t.A, Park Plaza, 59650
Villeneuve d'Ascq, France}

\address[authorlabel2]{Paul Painlevé  (CNRS, UMR 8524),  Université  de Lille 1, 59650, Villeneuve d'Ascq, France}

\address[authorlabel3]{Arts et Métiers
ParisTech centre de Lille, Laboratory of Applied Mathematics and Metrology (L2MA),  8 Boulevard Louis XIV, 59046 Lille Cedex,
France}
\address[authorlabel4]{LAGIS (CNRS, UMR 8146), \'{E}cole Centrale de Lille, BP 48,
Cit\'e Scientifique, 59650 Villeneuve d'Ascq, France}

\begin{abstract}
In this paper, the  numerical
differentiation by integration method based on Jacobi polynomials originally introduced by Mboup, Fliess and Join \cite{num,num0} is revisited in the central case where the used integration window is centered.  Such a method based on Jacobi polynomials was introduced through an algebraic approach \cite{num,num0} and extends the numerical
differentiation by integration method introduced by  Lanczos (1956) \cite{C. Lanczos}. The method  proposed here, rooted in \cite{num,num0}, is used to estimate the $n${th} ($n \in \mathbb{N}$) order derivative from
noisy data of a smooth function belonging to at least $C^{n+1+q}$
$(q \in \mathbb{N})$. In \cite{num,num0}, where the causal and anti-causal cases were investigated, the mismodelling due to the truncation of the Taylor expansion was investigated and improved allowing a small time-delay in the derivative estimation. Here, for the central case, we show that the bias error is $O(h^{q+2})$ where $h$ is the integration window length for $f\in C^{n+q+2}$ in the noise free case and the corresponding convergence
rate is $O(\delta^{\frac{q+1}{n+1+q}})$ where $\delta$ is the noise
level for a well-chosen integration window length.
Numerical examples show that this proposed method is stable and effective.

\end{abstract}


\end{frontmatter}

\section{Introduction}

Numerical differentiation is concerned with the numerical estimation
of derivatives of an unknown function (defined from $\mathbb{R}$ to $\mathbb{R}$) from its noisy measurement
data. It has attracted a lot of attention from different points of
view: observer design in the control literature \cite{R6,R19,R22},
digital filter in signal processing \cite{R5,R29}, the Volterra integral
equation of the first kind \cite{J. Cheng, R. Gorenflo} and
identification \cite{M. Hanke, Z. Wang2}. The problem of numerical
differentiation is ill-posed in the sense that a small error in
measurement data can induce a large error in the approximate
derivatives. Therefore, various numerical methods have been developed to obtain stable
algorithms more or less sensitive to additive noise. They mainly fall into five
categories: the finite difference methods \cite{I.R. Khan, R. Qu,
A.G. Ramm}, the mollification methods \cite{D.N. Hao, D.A. Murio,
D.A. Murio2}, the regularization methods \cite{G. Nakamura, T. Wei,
Y. Wang}, the algebraic methods \cite{num, num0} that are the roots of the results reported here and the
differentiation by integration methods \cite{C. Lanczos, S.K.
Rangarajana, Z. Wang}, i.e. using the Lanczos generalized
derivatives.

The Lanczos generalized derivative $D_{h}f$, defined in \cite{C.
Lanczos} by
\[
D_{h}f(x)=\frac{3}{2h^{3}}\int_{-h}^{h}t\,f(x+t)\,dt=\frac{3}{2h}\int
_{-1}^{1}t\,f(x+ht)\,dt,
\]
is an approximation to the first derivative of $f$ in the sense that
$D_{h}f(x)=f^{\prime}\left( x\right) + O(h^2) $. It is aptly called
a method of \emph{differentiation by integration}. Rangarajana and
al. \cite{S.K. Rangarajana} generalized it  for higher order
derivatives with
\[
D_{h}^{(n)}f(x)=\frac{1}{h^{n}}\int_{-1}^{1}\gamma_{n}L_{n}%
(t)\,f(x+ht)\,dt,\ n\in\mathbb{N},
\]
where $f$ is assumed to belong to $C^{n+2}(I)$ with $I$ being an
open interval of $\mathbb{R}$ and $L_{n}$ is the $n${th} order
Legendre polynomial. The coefficient $\gamma_{n}$ is equal to
$\frac{1\times3\times5\times \cdots\times(2n+1)}{2}$ and $2h>0$ is
the length of the integral window on which the estimates are
calculated. By applying the scalar product of the Taylor expansion
of $f$ at $x$\ with $L_{n}$  they showed that
$D_{h}^{(n)}f(x)=f^{(n)}(x)+O(h^{2})$. Recently, by using Richardson
extrapolation Wang and al. \cite{Z. Wang} have improved the
convergence rate for obtaining high order Lanczos derivatives with
the following affine schemes for any $n\in\mathbb{N}$
\[
D_{h,\lambda_{n}}^{(n)}f(x)=\frac{1}{h^{n}}\int_{-1}^{1}L_{n}(t)\left(
a_{n}\,f(x+ht)+b_{n}\,f(x+\lambda_{n}ht)\right)  dt,
\]
where $f$ is assumed to belong to $C^{n+4}(I)$, $a_{n}$, $b_{n}$ and
$\lambda_{n}$ are chosen such that $D_{h,\lambda_{n}}^{(n)}f(x)=f^{(n)}%
(x)+O(h^{4})$.

Very recently an algebraic setting for numerical differentiation of noisy signals was introduced in
\cite{FLI04a_compression_cras} and analyzed in \cite{num, num0}. The reader may find additional theoretical
foundations in \cite{ans, FLI03_Indent_ESAIM}. The
algebraic manipulations used in \cite{num, num0} are inspired by the one used in the
algebraic parametric estimation techniques \cite{mexico, mboup,
Med08}. Let us recall that \cite{num, num0} analyze a causal and an anti-causal version of numerical
differentiation by integration method based on Jacobi polynomials
\[
D_{h}^{(n)}f(x)=\frac{1}{(\pm h)^{n}}\int_{0}^{1}\gamma_{\kappa,\mu,n}\frac{d^n }{dt^n}
\left\{t^{\kappa+n}(1-t)^{\mu+n}\right\}\,f(x \pm ht)\,dt,\ n\in\mathbb{N},
\]
where $f$ is assumed to belong to $C^{n}(I)$ with $I$ being an open
interval of $\mathbb{R}$. The coefficient $\gamma_{\kappa,\mu,n}$ is
equal to $(-1)^n\frac{(\mu+\kappa+2n+1)!}{(\mu+n)!(\kappa+n)!}$,
where $\kappa,\mu$ are two integer parameters and $h>0$ is the
length of the integral window on which the estimates are calculated.
In \cite{num} the authors show that the mismodelling due to the
truncation of the Taylor expansion is improved allowing small
time-delay in the derivative estimation. Here in this article, we
propose to extend these \textit{differentiation by integration
}methods by using as in \cite{num, num0} Jacobi polynomials: for
this we use a central estimator (the integration window is now
$[-1,1]$) and the design parameters are now allowed  to be reals
which are strictly greater than $-1$. It is worth to mention that in
most of the practical applications the noise can be seen as an
integrable bounded function (which noise level is $\delta$ as  is
considered in this paper). Another point of view concerning the
noise definition/characterization is given in \cite{ans} for which
unbounded noise may appear. Let us mention that the Legendre
polynomials are one particular class of Jacobi polynomials, that
were used in \cite{S.K. Rangarajana} and \cite{Z. Wang} to obtain
higher order derivative estimations. Moreover, it can be seen that
these so obtained derivative estimators correspond to truncated
terms in the Jacobi orthogonal series. In fact, the choice of the
Jacobi polynomials comes from algebraic manipulations introduced in
the recent papers by M. Mboup, C. Join and M. Fliess \cite{num,
num0}, where the derivative estimations were  given by some
parameters in the causal and anti-causal cases. Here, we give the
derivative estimations in the central case with the same but
extended parameters used in \cite{num, num0}. If $f\in C^{n+q+2}$
then we show that the bias error is $O(h^{q+2})$ in the noise free
case (where $2h$ is the integration window length). We also show
that the corresponding convergence rate is
$O(\delta^{\frac{q+1}{n+1+q}})$ for a well-chosen integration window
length in the noisy case, where $\delta$ is the noise level. One can
see that the obtained causal estimators in \cite{num, num0} are well
suited for on-line estimation (which is of importance in signal
processing, automatic control, etc.) whereas here the proposed
central estimators are only suited for off-line applications. Let us
emphasize that those techniques exhibit good robustness properties
with respect to corrupting noises (see \cite{ans,shannon} for more
theoretical details). These robustness properties have already been
confirmed by numerous computer simulations and several laboratory
experiments. Hence, the robustness of the derivative estimators
presented in this paper can be ensured as shown by the 
results and simulations reported here.

This paper is organized as follows: in Section \ref{section2}
firstly a family of \textbf{central} estimators of the derivatives for higher orders
are introduced by using the $n${th} order Jacobi polynomials. The
corresponding convergence rate is $O(h)$ and can be improved to
$O(h^{2})$ when the Jacobi polynomials are {ultraspherical}
polynomials (see \cite{G. Szego}). Secondly, a new family of
estimators are given.\ They can be written as an affine combination
of the estimators proposed previously. Consequently, we show that if
$f\in C^{n+1+q}(I)$ with $q\in\mathbb{N}$ the corresponding
convergence rate is improved to $O(h^{q+1})$. Moreover, when the
Jacobi polynomials are  {ultraspherical} polynomials, if  $f\in
C^{n+2+q}(I)$ for any even integer $q$  the
corresponding convergence rate can be improved to $O(h^{q+2})$.
Numerical tests are given in Section \ref{section3} to verify the
efficiency and the stability of the proposed estimators.

\section{Derivative estimations by using Jacobi orthogonal series}

\label{section2}

Let $f^{\delta}=f+\varpi$ be a noisy function defined in an
open interval $I\subset\mathbb{R}$, where $f \in C^{n+1}(I)$ with
$n\in\mathbb{N}$ and the noise\footnote{More generally, the noise is
a stochastic process, which is bounded with certain probability and
integrable in the sense of convergence in mean square.} $\varpi$ is
bounded and integrable with a noise level $\delta$, $i.e.$
$\delta=\displaystyle\sup_{x\in I}|\varpi(x)|$. Contrary to
\cite{S.K. Rangarajana} where the $n${th} order Legendre
polynomials were used, we propose to use, as in \cite{num, num0}, the $n${th} order Jacobi
polynomial so as to obtain  estimates of the  $n${th} order
derivative of $f$. The $n${th} order Jacobi polynomials (see
\cite{G. Szego}) are defined as follows
\begin{equation} \label{jacobi}
P_{n}^{(\alpha,\beta)}(t)=\displaystyle\sum_{j=0}^{n}\binom{n+\alpha}{j}%
\binom{n+\beta}{n-j}\left(  \frac{t-1}{2}\right)  ^{n-j}\left(  \frac{t+1}%
{2}\right)  ^{j}%
\end{equation}
where $\alpha,\beta\in]-1,+\infty\lbrack$. Let us denote $\forall
g_{1},\, g_{2}\in C^{0}([-1,1])$, $\left\langle
g_{1}(\cdot),g_{2}(\cdot)\right\rangle
_{\alpha,\beta}=\int_{-1}^{1}w_{\alpha,\beta}(t)g_{1}(t)g_{2}(t)dt$,
where $w_{\alpha,\beta}(t)=(1-t)^{\alpha}(1+t)^{\beta}$ is the
weight function. Hence, we can denote its  associated norm  by
$\Vert\cdot\Vert_{\alpha,\beta}$.

We assume in this article that the parameter $h>0$\ and we denote
$I_{h}:=\left\{  x\in I;\left[  x-h,x+h\right]  \subset I\right\}
$.

\subsection{Minimal estimators}

In this subsection, let us ignore the noise $\varpi$ for a
moment. Then we can define a family of central estimators of $f^{(n)}$.

\begin{prop}
\label{proposition1} Let $f\in C^{n+1}(I)$, then a family of
central estimators of $f^{(n)}$ can be given as follows
\begin{equation}
\forall x\in I_{h},\ D_{h,\alpha,\beta}^{(n)}f(x)=\frac{1}{h^{n}}\int_{-1}%
^{1}\rho_{n,\alpha,\beta}(t)\,f(x+ht)dt,\label{d_h}%
\end{equation}
where $\rho_{n,\alpha,\beta}(t)=\frac{2^{-(n+\alpha+\beta+1)}n!}%
{B(n+\alpha+1,n+\beta+1)}P_{n}^{(\alpha,\beta)}(t)\,w_{\alpha,\beta}(t)$
with
$B(n+\alpha+1,n+\beta+1)=\frac{\Gamma(n+\alpha+1)\,\Gamma(n+\beta+1)}%
{\Gamma(2n+\alpha+\beta+2)}$.

Moreover, we have $D_{h,\alpha,\beta}^{(n)}f(x)=f^{(n)}(x)+O(h)$.
\end{prop}

\begin{remarque} \label{remarque1}
In order to compute $\rho_{n,\alpha,\beta}$, we should calculate
$P_{n}^{(\alpha,\beta)}$ whose computational complexity is
$O(n^2)$. Hence, the computational effort of $\rho_{n,\alpha,\beta}$
is $O(n^2)$.
\end{remarque}

\noindent\textbf{Proof.} By taking the Taylor expansion of $f$, we obtain for any $x \in I_h$ that there exists $\theta\in]x-h,x+h[$ such
that
\begin{equation} \label{Taylor}
f(x+ht)=f(x)+htf^{\prime}(x)+\cdots+\frac{h^{n}t^{n}}{n!}f^{(n)}%
(x)+\frac{h^{n+1}t^{n+1}}{(n+1)!}f^{(n+1)}(\theta).
\end{equation}
Substituting $(\ref{Taylor})$ in $(\ref{d_h})$, we deduce from the
classical orthogonal properties of the Jacobi polynomials (see
\cite{G. Szego}) that
\begin{align}
&  \int_{-1}^{1}\rho_{n,\alpha,\beta}(t)\,t^{m}\,dt=0,\ 0\leq m<n,\label{zero}%
\\
&  \int_{-1}^{1}\rho_{n,\alpha,\beta}(t)\,t^{n}\,dt=(n!).\label{n!}%
\end{align}
Using $(\ref{Taylor})$, $(\ref{zero})$ and $(\ref{n!})$, we can
conclude that
\[
D_{h,\alpha,\beta}^{(n)}f(x)=\frac{1}{h^{n}}\int_{-1}^{1}\rho_{n,\alpha,\beta
}(t)\,f(x+ht)dt=f^{(n)}(x)+O(h).
\]
Hence, this proof is completed.\hfill$\Box$\newline

In fact, we have taken an $n${th} order truncation in the Taylor
expansion of $f$ in Proposition \ref{proposition1} where $n$ is the
order of the estimated derivative. Thus, we call these estimators
 \emph{minimal estimators} (see \cite{num, num0}). Then, we can deduce the following
corollary.

\begin{corollaire}
Let $f \in C^{n+1}(I)$, then by assuming  that there exists
$M_{n+1}>0$ such that for any $x\in I$,
$\left|f^{(n+1)}(x)\right|\leq M_{n+1}$, we have
\begin{equation}
\left\Vert D_{h,\alpha,\beta}^{(n)}f(x)-f^{(n)}(x)\right\Vert
_{\infty}\leq C_{1}h,
\end{equation}
where
$C_{1}=\frac{M_{n+1}}{(n+1)!}\int_{-1}^{1}|t^{n+1}\rho_{n,\alpha,\beta
}(t)|\,dt$.
\end{corollaire}

When $\alpha=\beta$ the Jacobi polynomials are called ultraspherical
polynomials (see \cite{G. Szego}). In this case, we can improve the
convergence rate to $O(h^{2})$ by using the following lemma.

\begin{lemm}
\label{lemma} Let $P_{n}^{(\alpha,\alpha)}$  be the $n${th}
order ultraspherical polynomials, then we have
\begin{equation}
\label{zero2}\int_{-1}^{1} P_{n}^{(\alpha, \alpha)}(t)
w_{\alpha,\alpha}(t) \,t^{n+l} \,dt=0,
\end{equation}
where $l$ is an odd integer.
\end{lemm}

\noindent\textbf{Proof.} Recall the Rodrigues formula (see
\cite{G. Szego})
\begin{equation}\label{rodrigues}
P_{n}^{(\alpha,\beta)}(t)w_{\alpha,\beta}(t)=\frac{(-1)^{n}}{2^{n}n!}%
\frac{d^{n}}{dt^{n}}[w_{\alpha+n,\beta+n}(t)],
\end{equation}
we get, by substituting $(\ref{rodrigues})$ in $(\ref{zero2})$
and applying $n$ times integrations by parts, that
\begin{equation}\label{zero3}
\int_{-1}^{1}P_{n}^{(\alpha,\beta)}(t)w_{\alpha,\beta}(t)\,t^{n+l}%
\,dt=\frac{(n+l)!}{2^{n}(n!)^{2}}\int_{-1}^{1}w_{\alpha+n,\beta+n}%
(t)\,t^{l}\,dt.
\end{equation}
If $\alpha=\beta$ and $l$ is an odd number then $w_{\alpha+n,\beta
+n}(t)\,t^{l}$ is an odd function and the integral in
$(\ref{zero3})$ is equal to zero. Hence, this proof is completed.
\hfill$\Box$\newline

Consequently, we can deduce from Proposition \ref{proposition1} and
Lemma \ref{lemma} the following corollary.

\begin{corollaire}
\label{proposition11} Let $f\in C^{n+2}(I)$ and $\alpha=\beta$ in
Proposition \ref{proposition1}, then we obtain
\begin{equation}
\forall x\in I_{h},\
D_{h,\alpha,\alpha}^{(n)}f(x)=f^{(n)}(x)+O(h^{2}).
\end{equation}
Moreover, if we assume   that there exists ${M}_{n+2}>0$ such that
for any $x\in I$, $\left|f^{(n+2)}(x)\right|\leq {M}_{n+2}$,
then we have
\begin{equation}
\left\Vert D_{h,\alpha,\alpha}^{(n)}f(x)-f^{(n)}(x)\right\Vert _{\infty}%
\leq\hat{C}_{1}h^{2},
\end{equation}
where $\hat{C}_{1}=\frac{{M}_{n+2}}{(n+2)!}\int_{-1}^{1}|t^{n+2}%
\rho_{n,\alpha,\alpha}(t)|\,dt$.
\end{corollaire}

We can see in the following proposition  the relation between
minimal estimators of $f$ and  minimal estimators of $f^{(n)}$.

\begin{prop} \label{relation}
Let $f\in C^{n+1}(I)$, then we have
\begin{equation}
\forall x\in I_{h},\ D_{h,\alpha,\beta}^{(n)}f(x)=\frac{1}{(2h)^{n}}%
\frac{\Gamma(\alpha+\beta+2n+2)}{\Gamma(\alpha+\beta+n+2)}\sum_{j=0}%
^{n}(-1)^{n+j}\binom{n}{j}\ D_{h,\alpha_{n,j},\beta_{j}}^{(0)}f(x),
\end{equation}
where $\alpha_{n,j}=\alpha+n-j$ and $\beta_{j}=\beta+j$.
\end{prop}

In order to prove this proposition, we give the following lemma.
\begin{lemm}
\label{lemma2} For any $i\in\mathbb{N}$, we have
\begin{equation}
\forall x\in I_{h},\ \frac{\left\langle P_{i}^{(\alpha,\beta)}%
(t),f(x+ht)\right\rangle _{\alpha,\beta}}{\Vert P_{i}^{(\alpha,\beta)}%
\Vert^{2}_{\alpha,\beta}}=\sum_{j=0}^{i}(-1)^{i+j}\binom{i}{j}\frac{2i+\alpha+\beta
+1}{i+\alpha+\beta+1}\ D_{h,\alpha_{i,j},\beta_{j}}^{(0)}f(x),
\end{equation}
where $\alpha_{i,j}=\alpha+i-j$ and $\beta_{j}=\beta+j$.
\end{lemm}

\noindent\textbf{Proof.} Observe from the expression of the Jacobi
polynomials given in $(\ref{jacobi})$ that
\begin{equation}
\label{jacobi_w}P_{i}^{(\alpha,\beta)}(t)\,
w_{\alpha,\beta}(t)=\frac {1}{(-2)^{i}}\sum_{j=0}^{i}
\binom{i+\alpha}{j} \binom{i+\beta}{i-j} (-1)^{j}
w_{\alpha_{i,j},\beta_{j}}(t),
\end{equation}
we get
\begin{equation}%
\begin{split}
{\left\langle P_{i}^{(\alpha,\beta)}(t), f(x+ht) \right\rangle
_{\alpha,\beta }} = &  \frac{1}{(-2)^{i}}\sum_{j=0}^{i}
\binom{i+\alpha}{j} \binom{i+\beta }{i-j} (-1)^{j} \int_{-1}^{1}
w_{\alpha_{i,j},\beta_{j}}(t) \,f(x+ht) \,dt.
%= &  \frac{1}{(-2)^{i}}\sum_{j=0}^{i} \binom{i+\alpha}{j} \binom{i+\beta}{i-j}
%(-1)^{j}{\left\langle P_{0}^{(\alpha,\beta)}(t), f(x+ht) \right\rangle
%_{\alpha,\beta}}.
\end{split}
\end{equation}
Then, by using Proposition $\ref{proposition1}$ with $n=0$ and $P_{0}%
^{(\alpha_{i,j},\beta_{j})}(t)\equiv1$ we obtain
\begin{equation}%
\begin{split}
\label{pro}\frac{\left\langle P_{i}^{(\alpha,\beta)}(t), f(x+ht)
\right\rangle _{\alpha,\beta}}
{\|P_{i}^{(\alpha,\beta)}\|^{2}_{\alpha,\beta}} = &  \sum_{j=0}^{i}
\binom{i+\alpha}{j} \binom{i+\beta}{i-j} \frac{(-1)^{j}}{(-2)^{i}}
\frac{B(\alpha_{i,j}+1,\beta_{j}+1) 2^{\alpha_{i,j}+\beta_{j}+1}}%
{\|P_{i}^{(\alpha,\beta)}\|^{2}_{\alpha,\beta}}
D^{(0)}_{h,\alpha_{i,j},\beta_{j}}f(x).
\end{split}
\end{equation}
Recall that (see \cite{G. Szego})
\begin{equation} \label{norm}
\Vert P_{i}^{\left(  \alpha,\beta\right) }\Vert^{2}_{\alpha,\beta}=\frac{2^{\alpha+\beta+1}%
}{2i+\alpha+\beta+1}\frac{\Gamma(\alpha+i+1)\Gamma(\beta+i+1)}{\Gamma
(\alpha+\beta+i+1)\Gamma(i+1)},
\end{equation}
then  the proof is completed   by using $(\ref{norm})$ in
$(\ref{pro})$.
\hfill$\Box$\\

\noindent\textbf{Proof of Proposition \ref{relation}.}\ From
(\ref{d_h}), it is easy to show after some calculations
that%
\begin{equation} \label{relation2}
D_{h,\alpha,\beta}^{(n)}f(x)=\frac{1}{(2h)^{n}}\frac{\Gamma(\alpha
+\beta+2n+1)}{\Gamma(\alpha+\beta+n+1)}\frac{\left\langle
P_{n}^{(\alpha ,\beta)}(t),f(x+ht)\right\rangle
_{\alpha,\beta}}{\Vert P_{n}^{(\alpha,\beta
)}\Vert^{2}_{\alpha,\beta}}.
\end{equation}
Hence, this proof can be completed by using Lemma \ref{lemma2} and $(\ref{relation2})$. \hfill$\Box$\\

Now, we can see in the following proposition that the estimates given
in Proposition \ref{proposition1} are also equal to the first term
in the Jacobi orthogonal series expansion of $f^{(n)}(x+ht)$ at
point $t=0$.

\begin{prop}
\label{corollaire2} Let $f\in C^{n+1}(I)$, then the minimal
estimators of $f^{(n)}$ given in Proposition \ref{proposition1} can
be also written as follows
\begin{equation}
\forall x\in I_{h},\ D_{h,\alpha,\beta}^{(n)}f(x)=\frac{\left\langle
P_{0}^{(\alpha+n,\beta+n)}(t),f^{(n)}(x+ht)\right\rangle _{\alpha+n,\beta+n}%
}{\Vert P_{0}^{(\alpha+n,\beta+n)}\Vert^{2}_{\alpha+n,\beta+n}}\
P_{0}^{(\alpha+n,\beta+n)}(0).
\end{equation}
Moreover, we have
\begin{equation} \label{equalite}
\forall x\in I_{h},\ D_{h,\alpha,\beta}^{(n)}f(x)=D_{h,\alpha+n,\beta+n}%
^{(0)}f^{(n)}(x).
\end{equation}
\end{prop}

\noindent\textbf{Proof.} By using the Rodrigues formula in $(\ref{d_h})$
and applying $n$ times integrations by parts we get
\begin{equation*}
\begin{split}
D_{h,\alpha,\beta}^{(n)}f(x) & = \frac{1}{h^{n}} \frac{(-1)^{n}
2^{-(2n+\alpha +\beta+1)}}{B(n+\alpha+1,n+\beta+1)}\int_{-1}^{1}
\frac{d^{n}}{d t^{n}}[
w_{\alpha+n,\beta+n}(t)] \,f(x+ht) \,dt\\
& = \frac{
2^{-(2n+\alpha+\beta+1)}}{B(n+\alpha+1,n+\beta+1)}\int_{-1}^{1}
w_{\alpha+n,\beta+n}(t) \,f^{(n)}(x+ht) \,dt\\
& =D_{h,\alpha+n,\beta+n}^{(0)}f^{(n)}(x).
\end{split}
\end{equation*}
Then, by using $P_{0}^{(\alpha+n,\beta+n)}(t)\equiv1$ and $\|P_{0}%
^{(\alpha+n,\beta+n)}\|^{2}_{\alpha+n,\beta+n}=2^{2n+\alpha+\beta+1}{B(n+\alpha+1,n+\beta+1)}$,
we can achieve this proof.

\hfill$\Box$

%%%%%%%%%%%%%%%%%%%%%%%%%%%%%%%%%%%%%%%%%%%%%%%%%%%%%%%%%%%%%%%%%%%%%%%%%%%%%%%%%%%%%%%%%%%%%%%%%%%%%%%%%%%%%%%%%%%%%%%%%%%
\subsection{Affine estimators}

It is shown in Proposition \ref{corollaire2} that the minimal
estimators of $f^{(n)}(x)$ given in Proposition \ref{proposition1}
are equal to the value of the $0$ order truncated Jacobi orthogonal
series expansion of $f^{(n)}(x+ht)$ at $t=0$. Let us assume that
$f\in C^{n+1}(I)$, then we define now the $q${th} ($q \in
\mathbb{N}$) order truncated Jacobi orthogonal series of
$f^{(n)}(x+ht)$ by the following operator
\begin{equation} \label{orthogonal_series}
\forall x\in I_{h},\ D_{h,\alpha,\beta,q}^{(n)}f(x+th):=\sum_{i=0}^{q}%
\frac{\left\langle
P_{i}^{(\alpha+n,\beta+n)}(\cdot),f^{(n)}(x+h \cdot)\right\rangle
_{\alpha+n,\beta+n}}{\Vert P_{i}^{(\alpha+n,\beta+n)}\Vert^{2}_{\alpha+n,\beta+n}}\ P_{i}%
^{(\alpha+n,\beta+n)}(t).
\end{equation}

Take $t=0$ in $(\ref{orthogonal_series})$, we obtain a family of
estimators of $f^{(n)}(x)$ with
\begin{equation} \label{d_h_N}
\forall x\in I_{h},\ D_{h,\alpha,\beta,q}^{(n)}f(x)=\sum_{i=0}^{q}%
\frac{\left\langle
P_{i}^{(\alpha+n,\beta+n)}(\cdot),f^{(n)}(x+h \cdot)\right\rangle
_{\alpha+n,\beta+n}}{\Vert P_{i}^{(\alpha+n,\beta+n)}\Vert^{2}_{\alpha+n,\beta+n}}\ P_{i}%
^{(\alpha+n,\beta+n)}(0).
\end{equation}

To better explain  our method, let us recall some well-known facts.
We consider the subspace of $C^{0}([-1,1])$, defined by
\begin{equation}\label{}
    \mathcal{H}_q=\text{span}\left\{P_{0}^{(\alpha+n,\beta+n)},P_{1}^{(\alpha+n,\beta+n)},\cdots,P_{q}^{(\alpha+n,\beta+n)}\right\}.
\end{equation}
Equipped with the inner product $\left\langle
\cdot,\cdot\right\rangle_{\alpha+n,\beta+n}$, $\mathcal{H}_q$ is clearly a reproducing kernel Hilbert space  \cite{Aronszajn}, \cite{Saitoh}, with the reproducing kernel
\begin{equation}\label{}
    \mathcal{K}_q(\tau,t)=\sum_{i=0}^{q}%
\frac{
P_{i}^{(\alpha+n,\beta+n)}(\tau) P_{i}%
^{(\alpha+n,\beta+n)}(t)}{\Vert P_{i}^{(\alpha+n,\beta+n)}\Vert^{2}_{\alpha+n,\beta+n}}.
\end{equation}
The reproducing property implies that for any function $f^{(n)}(x+h \cdot)$ belonging to $C^{0}([-1,1])$, we have
\begin{equation}\label{}
    \left\langle
\mathcal{K}_q(\cdot,t),f^{(n)}(x+h \cdot)\right\rangle_{\alpha+n,\beta+n}=D_{h,\alpha,\beta,q}^{(n)}f(x+th),
\end{equation}
where $D_{h,\alpha,\beta,q}^{(n)}f(x+h \cdot)$ stands for the orthogonal projection of $f^{(n)}(x+h \cdot)$ on $\mathcal{H}_q$. Thus, the estimators given in $(\ref{d_h_N})$
can be obtained by taking $t=0$.

We will see in the following proposition that the estimators
$D_{h,\alpha ,\beta,q}^{(n)}f(x)$ can be written as an affine
combination of different minimal estimators. These estimators are
called \emph{affine estimators} as in \cite{num}.

\begin{prop}
\label{proposition5} Let $f\in C^{n+1}(I)$,
then we have
\begin{equation} \label{affine_relation}
\forall x\in I_{h},\ D_{h,\alpha,\beta,q}^{(n)}f(x)=\sum_{i=0}^{q}%
{P_{i}^{(\alpha+n,\beta+n)}(0)}\frac{2i+\alpha+\beta+2n+1}{i+\alpha+\beta+2n+1} \sum_{j=0}^{i}(-1)^{i+j}\binom{i}{j}%
 D_{h,\alpha_{i,j},\beta_{j}%
}^{(n)}f(x),
\end{equation}
where  $q\in\mathbb{N}$, $\alpha_{i,j}=\alpha+i-j$ and $\beta_{j}=\beta+j$. Moreover, we have
\begin{equation} \label{affine_sum}
\sum_{i=0}^{q}%
{P_{i}^{(\alpha+n,\beta+n)}(0)}\frac{2i+\alpha+\beta+2n+1}{i+\alpha+\beta+2n+1} \sum_{j=0}^{i}(-1)^{i+j}\binom{i}{j}%
=1.
\end{equation}
\end{prop}

\noindent\textbf{Proof.} By replacing $\alpha$ by $\alpha+n$, $\beta$
by $\beta+n$ and $f(x+ht)$ by $f^{(n)}(x+ht)$  in Lemma
\ref{lemma2}, we obtain
\begin{equation}
\frac{\left\langle P_{i}^{(\alpha+n,\beta+n)}(t), f^{(n)}(x+ht)
\right\rangle _{\alpha+n,\beta+n}}
{\|P_{i}^{(\alpha+n,\beta+n)}\|^{2} _{\alpha
+n,\alpha+n}}=\sum_{j=0}%
^{i}(-1)^{i+j} \binom{i}{j}
\frac{2i+\alpha+\beta+2n+1}{i+\alpha+\beta+2n+1} \
D_{h,\alpha_{i,j}+n,\beta_{j}+n}^{(0)}f^{(n)}(x).
\end{equation}
Then $(\ref{affine_relation})$ can be obtained  by using $(\ref{equalite})$ and
$(\ref{d_h_N})$. By using the Binomial relation, $(\ref{affine_sum})$ can be easily obtained. \hfill$\Box$\newline

Hence, by using Proposition \ref{proposition1} an explicit
formulation of these affine estimators is obtained in the following
corollary.
\begin{corollaire}
\label{proposition2} Let $f\in C^{n+1}(I)$,
then the affine estimators of $f^{(n)}$ can be written as
\begin{equation} \label{D_h_N}
\forall x\in I_{h},\
D_{h,\alpha,\beta,q}^{(n)}f(x)=\frac{1}{h^{n}}\int
_{-1}^{1}Q_{\alpha,\beta,n,q}(t)\,f(x+ht)\,dt,
\end{equation}
where
\begin{equation} \label{Q}
Q_{\alpha,\beta,n,q}(t)=\displaystyle\sum_{i=0}^{q}{P_{i}^{(\alpha+n,\beta
+n)}(0)}\sum_{j=0}^{i}(-1)^{i+j}\binom{i}{j}\frac{2i+\alpha+\beta
+2n+1}{i+\alpha+\beta+2n+1}\,\rho_{n,\alpha_{i,j},\beta_{j}}(t)
\end{equation}
with $q\in\mathbb{N}$, $\rho_{n,\alpha_{i,j},\beta_{j}}$ given in Proposition
\ref{proposition1} and $\alpha_{i,j}=\alpha+i-j$,
$\beta_{j}=\beta+j$.
\end{corollaire}

Consequently these affine estimators are also
\textit{differentiation by integration} estimators.

\begin{remarque} \label{remarque2}
$Q_{\alpha,\beta,n,q}$ is a sum of  $\frac{1}{2}(q+1)(q+2)$ terms.
According to Remark \ref{remarque1}, the computational effort of
each term is  $O(n^2)$. Hence, the computational effort of
$Q_{\alpha,\beta,n,q}$ is also $O(n^2)$.
\end{remarque}

It is shown in Proposition \ref{proposition1} that the convergence
rate of minimal estimators is $O(h)$. We will see in the following
proposition that the convergence rate of affine estimators can be
improved to $O(h^{q+1})$.

\begin{prop}
\label{proposition3} Let $f\in C^{n+1+q}(I)$ with $q\in\mathbb{N}$,
then we have
\begin{equation}
\forall x\in I_{h},\
D_{h,\alpha,\beta,q}^{(n)}f(x)=f^{(n)}(x)+O(h^{q+1}).
\end{equation}
Moreover, if we assume that there exists $M_{n+1+q}>0$ such that for
any $x\in I$, $\left|f^{(n+1+q)}(x)\right|\leq M_{n+1+q}$, then
we have
\begin{equation}
\left\Vert D_{h,\alpha,\beta,q}^{(n)}f(x)-f^{(n)}(x)\right\Vert
_{\infty }\leq C_{2}h^{q+1},
\end{equation}
where
$C_{2}=\frac{M_{n+1+q}}{(n+1+q)!}\int_{-1}^{1}|t^{n+1+q}Q_{\alpha,\beta
,n,q}(t)|\,dt$.
\end{prop}

\noindent\textbf{Proof.} By taking the Taylor expansion of $f$, we
get for any $x \in I_h$ there exists $\xi\in]x-h,x+h[$ such that
\begin{equation}
\label{Taylor2}f(x+ht)= f_{n+q}(x+ht) + \frac{h^{n+1+q}
t^{n+1+q}}{(n+1+q)!} f^{(n+1+q)}(\xi),
\end{equation}
where $f_{n+q}(x+ht)=\displaystyle\sum_{j=0}^{n+q} \frac{h^{j}
t^{j}}{j!} f^{(j)}(x)$ is the $(n+q)${th} order truncated Taylor
series expansion of $f(x+ht)$.

Let us take the Jacobi orthogonal series expansion of
$f_{n+q}^{(n)}(x+ht)$. Then by taking $t=0$, we obtain
\begin{equation}%
\begin{split} \label{f_q}
f_{n+q}^{(n)}(x) = & \sum_{i=0}^{q}\frac{\left\langle P_{i}^{(\alpha
+n,\beta+n)}(t), f_{n+q}^{(n)}(x+ht) \right\rangle
_{\alpha+n,\beta+n}} {\|P_{i}^{(\alpha+n,\beta+n)}\|^{2} _{\alpha
+n,\alpha+n}} \ P_{i}^{(\alpha+n,\beta+n)}(0).
\end{split}
\end{equation}

Similarly  to  $(\ref{d_h_N})$  we obtain
\begin{equation}
\label{f_N}f_{n+q}^{(n)}(x)=   \frac{1}{h^{n}} \int_{-1}^{1}
Q_{\alpha ,\beta,n,{q}}(t) f_{n+q}(x+ht) dt,
\end{equation}
from $(\ref{D_h_N})$ where $Q_{\alpha,\beta,n,q}$ is given in
Corollary \ref{proposition2} by $(\ref{Q})$.

By calculating the value of the $n${th} order derivative of  $f_{n+q}^{(n)}$ at $t=0$, we obtain $f_{n+q}^{(n)}(x)=f^{(n)}(x)$. Then by using
$(\ref{D_h_N})$ and $(\ref{f_N})$ we obtain
\begin{equation*}
\begin{split}
D_{h,\alpha,\beta,{q}}^{(n)}f(x) - f^{(n)}(x) = &  \frac{1}{h^{n}}
\int
_{-1}^{1} Q_{\alpha,\beta,n,q}(t) \left[ f(x+ht)-f_{n+q}(x+ht) \right]  dt\\
= & \frac{ h^{q+1}}{(n+1+q)!} \int_{-1}^{1} Q_{\alpha,\beta,n,q}(t)
t^{n+1+q}
f^{(n+1+q)}(\xi) dt\\
= &  O(h^{q+1}).
\end{split}
\end{equation*}

Consequently, if for any $x \in I$ $\left|f^{(n+1+q)}(x)
\right|\leq M_{n+1+q}$, then we have
\begin{equation*}
\left\|  D_{h,\alpha,\beta,{q}}^{(n)}f(x) - f^{(n)}(x)\right\|
_{\infty} \leq h^{q+1} \frac{M_{n+1+q}}{(n+1+q)!}\int_{-1}^{1}
|t^{n+1+q} Q_{\alpha,\beta,n,q}(t)| \, dt.
\end{equation*}
\hfill$\Box$

%\begin{remarque}
We can deduce that the affine estimator for $f^{(n)}(x)$ obtained by
taking the $q${th} order truncated Jacobi orthogonal series
expansion of $f^{(n)}(x+h\cdot)$ can be also obtained by taking the
$(n+q)${th} order truncated  Taylor series expansion of $f$ with a
scalar product of Jacobi polynomials.

Moreover, let $f_{n+q}(x+ht)=f_{n}(x+ht)+r_{q}(x+ht)$ where
$r_{q}(x+ht)=\displaystyle\sum_{j=n+1}^{n+q} \frac{h^{j} t^{j}}{j!}
f^{(j)}(x)$ for $q \geq 1$ and $r_{q}(x+ht)=0$ for $q=0$, then
$(\ref{f_q})$ becomes
\begin{equation*}%
f_{n+q}^{(n)}(x) = \sum_{i=0}^{q}\frac{\left\langle P_{i}^{(\alpha
+n,\beta+n)}(t), f_{n}^{(n)}(x+ht) \right\rangle
_{\alpha+n,\beta+n}} {\|P_{i}^{(\alpha+n,\beta+n)}\|^{2} _{\alpha
+n,\alpha+n}} \ P_{i}^{(\alpha+n,\beta+n)}(0)+ R,
\end{equation*}
where $R=\displaystyle\sum_{i=0}^{q}\frac{\left\langle
P_{i}^{(\alpha +n,\beta+n)}(t), r_{q}^{(n)}(x+ht) \right\rangle
_{\alpha+n,\beta+n}} {\|P_{i}^{(\alpha+n,\beta+n)}\|^{2} _{\alpha
+n,\alpha+n}} \ P_{i}^{(\alpha+n,\beta+n)}(0)$.

Observe that $f_{n}^{(n)}(x+h\cdot)$ is a $0${th} order polynomial, then by using the orthogonal properties  of
$P_{i}^{(\alpha +n,\beta+n)}$ we have
\begin{equation*}%
\begin{split}
&\sum_{i=0}^{q}\frac{\left\langle P_{i}^{(\alpha +n,\beta+n)}(t),
f_{n}^{(n)}(x+ht) \right\rangle _{\alpha+n,\beta+n}}
{\|P_{i}^{(\alpha+n,\beta+n)}\|^{2} _{\alpha +n,\alpha+n}} \
P_{i}^{(\alpha+n,\beta+n)}(0) \\
=& \frac{\left\langle P_{0}^{(\alpha +n,\beta+n)}(t),
f_{n}^{(n)}(x+ht) \right\rangle _{\alpha+n,\beta+n}}
{\|P_{0}^{(\alpha+n,\beta+n)}\|^{2} _{\alpha +n,\alpha+n}} \
P_{0}^{(\alpha+n,\beta+n)}(0)=f_{n}^{(n)}(x).
\end{split}
\end{equation*}

By calculating the value of the $n${th} order derivative of  $f_{n+q}^{(n)}$  and $f_{n}^{(n)}$ at $t=0$, we obtain $f_{n+q}^{(n)}(x)=f_{n}^{(n)}(x)=f^{(n)}(x)$. Hence, we get  $R=0$.
Hence, we can deduce that
\begin{equation} \label{r}
R=\frac{1}{h^{n}} \int_{-1}^{1} Q_{\alpha ,\beta,n,{q}}(t)
r_{q}(x+ht) dt=0,
\end{equation}
where $Q_{\alpha,\beta,n,q}$ is given in Corollary
\ref{proposition2} by $(\ref{Q})$.

Consequently, $(\ref{r})$ explains why the convergence rate can be
improved from $O(h)$ to $O(h^{q+1})$: the price to pay is some more
smoothness hypotheses  on the function $f$.
%\end{remarque}

If we consider the noisy function $f^{\delta}$, then it is
sufficient to replace $f(x+h\cdot)$ in $(\ref{D_h_N})$ by
$f^{\delta}(x+h\cdot)$ so as to estimate $f^{(n)}(x)$. Then we have the
following definition.

\begin{definition}
\label{proposition4} Let $f^{\delta}=f+\varpi$ be a noisy
function, where $f \in C^{n+1}(I)$ and $\varpi$ is a bounded and
integrable noise with a noise level $\delta$. Then a family of
estimators of $f^{(n)}$ is defined as follows
\begin{equation}
\forall x \in I_h, \ D_{h,\alpha,\beta,{q}}^{(n)}f^{\delta}(x)=\frac{1}{h^{n}}\int_{-1}%
^{1}Q_{\alpha,\beta,n,{q}}(t)\,f^{\delta}(x+ht)\,dt,
\end{equation}
where $Q_{\alpha,\beta,n,q}$ is given by $(\ref{Q})$.
\end{definition}

In the following proposition we study the estimation error for these estimators.
\begin{prop}
\label{corollary3} Let $f^{\delta}$ be a noisy function where $f
\in C^{n+1+q}(I)$ and $\varpi$ is a bounded and integrable noise
with a noise level $\delta$, then
\begin{equation}
\left\Vert
D_{h,\alpha,\beta,{q}}^{(n)}f^{\delta}(x)-f^{(n)}(x)\right\Vert
_{\infty}\leq C_{2}h^{q+1}+C_{3}\frac{\delta}{h^{n}},
\end{equation}
where $C_{2}$ is given in Proposition \ref{proposition3} and
$C_{3}=\int _{-1}^{1}|Q_{\alpha,\beta,n,q}(t)|\,dt.$

Moreover, if we choose $h=\left[
\frac{nC_{3}}{(q+1)C_{2}}\delta\right] ^{\frac{1}{n+q+1}}$, then we
have
\begin{equation}
\left\Vert
D_{h,\alpha,\beta,{q}}^{(n)}f^{\delta}(x)-f^{(n)}(x)\right\Vert
_{\infty}=O(\delta^{\frac{q+1}{n+1+q}}).
\end{equation}
\end{prop}

\noindent\textbf{Proof.} Since
\begin{equation*}
\left\|  D_{h,\alpha,\beta,{q}}^{(n)}f^{\delta}(x) -
D_{h,\alpha,\beta
,{q}}^{(n)}f(x)\right\| _{\infty} =\left\|  D_{h,\alpha,\beta,{q}}%
^{(n)}\left[ f^{\delta}(x) - f(x)\right] \right\|
_{\infty}\leq\frac{\delta }{h^{n}} \int_{-1}^{1} \left|
Q_{\alpha,\beta,n,q}(t)\right|  \, dt,
\end{equation*}
by using Proposition \ref{proposition3} we get
\begin{equation*}
\begin{split}
\left\|  D_{h,\alpha,\beta,{q}}^{(n)}f^{\delta}(x) -
f^{(n)}(x)\right\| _{\infty}  & \leq\left\|
D_{h,\alpha,\beta,{q}}^{(n)}f^{\delta}(x) -
D_{h,\alpha,\beta,{q}}^{(n)}f(x)\right\| _{\infty}+\left\|
D_{h,\alpha
,\beta,{q}}^{(n)}f(x) - f^{(n)}(x)\right\| _{\infty}\\
& \leq C_{2} h^{q+1} + C_{3} \frac{\delta}{h^{n}},
\end{split}
\end{equation*}
where $C_{3}=\int_{-1}^{1} | Q_{\alpha,\beta,n,q}(t)| \, dt.$  Let
us denote the
error bound by $\psi(h)=C_{2}h^{q+1}+C_{3}\frac{\delta}{h^{n}}$. Consequently, we can calculate its minimum value.  It
is obtained for $h^{\ast}=\left[
\frac{nC_{3}}{(q+1)C_{2}}\delta\right]^{\frac {1}{n+q+1}}$ and
\begin{equation}
\psi(h^{\ast})=\frac{n+1+q}{q+1}\left(  \frac{q+1}{n}\right) ^{\frac
{n}{n+1+q}}C_{2}^{\frac{n}{n+1+q}}C_{3}^{\frac{q+1}{n+1+q}}\delta^{\frac
{q+1}{n+1+q}}.
\end{equation}
Then, the proof is completed. \hfill$\Box$\\

%\begin{remarque}
In Proposition \ref{proposition5},  we improve the convergence rate
 from $O(h)$ to $O(h^{q+1})$ ($q \in
\mathbb{N}$) for the exact function $f$ by taking an affine
combination of  minimal estimators of $f^{(n)}$. Here, the
convergence rate is also improved for  noisy functions. It passes
from $O(\delta^{\frac{1}{n+1}})$ to $O(\delta^{\frac{q+1}{n+1+q}})$
if we choose $h=c \,\delta^{\frac{1}{n+1+q}}$, where $c$ is a
constant.
%\end{remarque}

\begin{remarque}
Usually, the sampling data are given in discrete case. We
should use a numerical integration method to approximate the
integrals in our estimators. This numerical method will produce a
numerical error. Hence, we  always set the value of $h$ larger than
the optimal one calculated in the previous proof.
\end{remarque}

We have seen in the previous subsection that the convergence rate of
 minimal estimators can be improved to $O(h^{2})$ when
$\alpha=\beta$. Let us then study the convergence rate of affine
estimators in this case.

\begin{corollaire}
\label{corollary4} Let $f \in C^{n+2+q}(I)$ where $q$ is an even
integer. If we set  $\alpha=\beta$ in
(\ref{d_h_N}), then we have
\begin{equation}
\forall x\in I_{h},\
D_{h,\alpha,\alpha,q}^{(n)}f(x)=f^{(n)}(x)+O(h^{q+2}).
\end{equation}
Moreover, if we assume that there exists ${M}_{n+2+q}>0$ such that
for any $x\in I$, $\left|f^{(n+q+2)}(x)\right|\leq {M}_{n+2+q}$,
then we have
\begin{equation} \label{bound}
\left\Vert D_{h,\alpha,\alpha,q}^{(n)}f(x)-f^{(n)}(x)\right\Vert
_{\infty }\leq\hat{C}_{2}h^{q+2},
\end{equation}
where $\hat{C}_{2}=\frac{{M}_{n+2+q}}{(n+q+2)!}\int_{-1}^{1}|t^{n+q+2}%
Q_{\alpha,n,q}(t)|\,dt$ and
\begin{equation} \label{Q2}
Q_{\alpha,n,q}(t)=\displaystyle\sum_{i=0}^{\frac{q}{2}}{P_{2i}^{(\alpha+n,\alpha
+n)}(0)}\sum_{j=0}^{2i}(-1)^{j}\binom{2i}{j}\frac{4i+2\alpha
+2n+1}{2i+2\alpha+2n+1}\,\rho_{n,\alpha_{2i,j},\beta_{j}}(t)
\end{equation}
with $\rho_{n,\alpha_{2i,j},\beta_{j}}$ given in Proposition
\ref{proposition1} and $\alpha_{2i,j}=\alpha+2i-j$,
$\beta_{j}=\alpha+j$.
\end{corollaire}

\noindent\textbf{Proof.} Observe that $P_{q+1}^{(\alpha+n, \alpha
+n)}(-t)=(-1)^{(q+1)}P_{q+1}^{(\alpha+n, \alpha+n)}(t)$ for any $t \in [-1,1]$ (see \cite{G.
Szego} p.80), we obtain $P_{q+1}^{(\alpha+n, \alpha+n)}(0)=0$. Hence,
$(\ref{d_h_N})$ becomes
\begin{equation*}
\begin{split}
D_{h,\alpha,\alpha,q}^{(n)}f(x)  &
=\sum_{i=0}^{q+1}\frac{\left\langle P_{i}^{(\alpha+n, \alpha+n)}(t),
f^{(n)}(x+ht) \right\rangle _{\alpha +n,\alpha+n}}
{\|P_{i}^{(\alpha+n,\alpha+n)}\|^{2} _{\alpha +n,\alpha+n}} \
P_{i}^{(\alpha +n,\alpha+n)}(0).
\end{split}
\end{equation*}

If $f \in C^{n+2+q}(I)$, then let us take $f_{n+q+1}$ as the
$(n+q+1)${th} order truncated Taylor series expansion of $f(x+ht)$.
By taking the Jacobi orthogonal series expansion of
$f^{(n)}_{n+q+1}$
\begin{equation*}
f^{(n)}_{n+q+1}(x)   =\sum_{i=0}^{q+1}\frac{\left\langle
P_{i}^{(\alpha+n, \alpha+n)}(t), f^{(n)}_{n+q+1}(x+ht) \right\rangle
_{\alpha+n,\alpha+n}} {\|P_{i}^{(\alpha+n,\alpha+n)}\|^{2} _{\alpha
+n,\alpha+n}} \ P_{i}^{(\alpha+n,\alpha+n)}(0),
\end{equation*}
we obtain
\begin{equation*}
\begin{split}
D_{h,\alpha,\alpha,q}^{(n)}f(x)-f^{(n)}(x)  & = D_{h,\alpha,\alpha
,q}^{(n)}f(x)-f^{(n)}_{n+q+1}(x)\\
& =\sum_{i=0}^{q+1}\frac{\left\langle P_{i}^{(\alpha+n,
\alpha+n)}(t), f^{(n)}(x+ht)-f^{(n)}_{n+q+1}(x+ht) \right\rangle
_{\alpha+n,\alpha+n}} {\|P_{i}^{(\alpha+n,\alpha+n)}\|^{2} _{\alpha
+n,\alpha+n}} \ P_{i}^{(\alpha+n,\alpha+n)}(0)\\
& =\frac{1}{h^{n}} \int_{-1}^{1} Q_{\alpha,\alpha,n,q}(t) \left[
f(x+ht)-f_{n+q+1}(x+ht) \right]  dt\\
& =\frac{ h^{q+2}}{(n+2+q)!} \int_{-1}^{1} Q_{\alpha,\alpha,n,q}(t)
t^{n+2+q}
f^{(n+2+q)}(\xi^{\prime}) dt, \ \ \xi^{\prime}\in]x-h, x+h[\\
& = O(h^{q+2}).
\end{split}
\end{equation*}
Consequently, $(\ref{bound})$ follows directly from the {hypothesis}
on $\left| f^{(n+q+2)}(x) \right|$. Since $P_{i}^{(\alpha+n,
\alpha+n)}(0)=0$ for any odd integer $i$, $(\ref{Q2})$ can be
obtained by using $(\ref{Q})$. Then this proof is completed.
\hfill$\Box$

\begin{remarque} \label{remarque3}
According to \cite{Poffald}, we can deduce the asymptotic behavior
 of the number $\xi'$ when $h \rightarrow 0^+$
\begin{equation}\label{theta}
    \lim_{h \rightarrow 0^+} \frac{|\xi'-x|}{h}=\frac{1}{n+q+3}.
\end{equation}
\end{remarque}

Similarly to Proposition  \ref{corollary3}, we can obtain the
following corollary.
\begin{corollaire} \label{corollary5}
Let $f\in C^{n+2+q}(I)$ where $q$ is an even integer. If
$\alpha=\beta$ in Definition  \ref{proposition4}, then the
estimation error for $D_{h,\alpha,\alpha,{q}}^{(n)}f^{\delta}(x)$ is
given by
\[
\left\Vert
D_{h,\alpha,\alpha,{q}}^{(n)}f^{\delta}(x)-f^{(n)}(x)\right\Vert
_{\infty}\leq\hat{C}_{2}h^{q+2}+C_{3}\frac{\delta}{h^{n}},
\]
where $\hat{C}_{2}$ is given in Corollary \ref{corollary4} and
$C_{3}$ is given in Proposition \ref{corollary3}.

Moreover, if we choose $\hat{h}=\left[
\frac{nC_{3}}{(q+2)\hat{C}_{2}}\delta\right] ^{\frac{1}{n+q+2}}$,
then we have
\[
\left\Vert
D_{h,\alpha,\alpha,{q}}^{(n)}f^{\delta}(x)-f^{(n)}(x)\right\Vert
_{\infty}=O(\delta^{\frac{q+2}{n+2+q}}).
\]
\end{corollaire}

In the following proposition, if we assume that $f \in C^{n-1}(I)$
then we can define the generalized derivative of $f^{(n)}$. We can
see that if the right and left hand derivatives for the $n${th}
order exist, then this generalized derivative converges to the
average value of these one-sided derivatives.

\begin{prop}\label{proposition6}
Let $f \in C^{n-1}(I)$, then we define the generalized derivative of
$f^{(n)}$ by
\begin{equation} \label{D_h_N7}
\forall x\in I_{h},\
D_{h,\alpha,\alpha,q}^{(n)}f(x)=\frac{1}{h^{n}}\int
_{-1}^{1}Q_{\alpha,n,q}(t)\,f(x+ht)\,dt,
\end{equation}
where $Q_{\alpha,n,q}$ is defined by (\ref{Q2}). Moreover,
 if  $f^{(n)}_+(x)$ and
$f^{(n)}_-(x)$ exist at any point $x \in I_h$, then we have
\begin{equation}\label{}
    \lim_{h \rightarrow 0^+}
D_{h,\alpha,\alpha,{q}}^{(n)}f(x)=\frac{1}{2}\left(f^{(n)}_+(x)+f^{(n)}_-(x)\right),
\end{equation}
where $f^{(n)}_+$ (resp. $f^{(n)}_-$) denotes the right (resp. left)
hand  derivative for the $n${th} order.
\end{prop}

Before proving this proposition, let us give the following lemma.

\begin{lemm} \label{lemma3}
Let $n \in \mathbb{N}$ and $Q_{\alpha,n,q}$ be the function defined
on $[-1,1]$ by $(\ref{Q2})$  where $q$ is an even integer.
If $n$ is  even then $Q_{\alpha,n,q}$ is also even, odd else.
\end{lemm}

\noindent\textbf{Proof.}
By taking $\alpha=\beta$ in (\ref{d_h_N}), we obtain
\begin{equation} \label{d_h_N2}
\forall x\in I_{h},\ D_{h,\alpha,\alpha,q}^{(n)}f(x)=\sum_{i=0}^{q}%
\frac{P_{i}%
^{(\alpha+n,\alpha+n)}(0)}{\Vert P_{i}^{(\alpha+n,\alpha+n)}\Vert^{2}_{\alpha+n,\alpha+n}}\int_{-1}^1
P_{i}^{(\alpha+n,\alpha+n)}(t) w_{\alpha+n,\alpha+n}(t) f^{(n)}(x+h t) dt.
\end{equation}
By using (\ref{jacobi_w}) and replacing $\alpha,\beta$ by $\alpha+n$, we get for $l=0,\cdots,n-1$
\begin{equation*}
\frac{d^l}{d t^l}\left[P_{i}^{(\alpha+n,\alpha+n)}(t)\,
w_{\alpha+n,\alpha+n}(t)\right]=\frac {1}{(-2)^{i}}\sum_{j=0}^{i}
\binom{i+\alpha+n}{j} \binom{i+\alpha+n}{i-j} (-1)^{j}
\frac{d^l}{d t^l} \left[w_{\alpha_{i+n,j},\alpha_{j+n}}(t)\right],
\end{equation*}
where $\alpha_{i+n,j}=\alpha+i+n-j$,
$\alpha_{j+n}=\alpha+j+n$.
Then, by applying the Rodrigues formula, we get
\begin{equation*}%\label{valeur}
\begin{split}
&\frac{d^l}{d t^l}\left[P_{i}^{(\alpha+n,\alpha+n)}(t)\,
w_{\alpha+n,\alpha+n}(t)\right]\\=&\frac {l!}{2^{i-l}}\sum_{j=0}^{i}
\binom{i+\alpha+n}{j} \binom{i+\alpha+n}{i-j} (-1)^{j+i+l}
P_{n}^{(\alpha_{i+n-l,j},\alpha_{j+n-l})}(t)w_{\alpha_{i+n-l,j},\alpha_{j+n-l}}(t),
\end{split}
\end{equation*}
where $\alpha_{i+n-l,j}=\alpha+i+n-j-l$,
$\alpha_{j+n-l}=\alpha+j+n-l$.  Hence,
we get that $\frac{d^l}{d t^l}[P_{i}^{(\alpha+n,\alpha+n)}\,
w_{\alpha+n,\alpha+n}]$ are equal to $0$ at $-1$ and $1$. Thus, by applying $n$ times integrations by parts in $(\ref{d_h_N2})$,
we obtain
\begin{equation} \label{d_h_N3}
\forall x\in I_{h},\ D_{h,\alpha,\alpha,q}^{(n)}f(x)=\frac{(-1)^n}{h^n}\sum_{i=0}^{q}%
\frac{P_{i}%
^{(\alpha+n,\alpha+n)}(0)}{\Vert P_{i}^{(\alpha+n,\alpha+n)}\Vert^{2}_{\alpha+n,\alpha+n}}\int_{-1}^1
\frac{d^n}{d t^n} \left[ P_{i}^{(\alpha+n,\alpha+n)}(t) w_{\alpha+n,\alpha+n}(t)\right] f(x+h t) dt.
\end{equation}
By using Corollary \ref{proposition2} with $\alpha=\beta$, we get
\begin{equation}\label{Q3}
    Q_{\alpha,n,q}(t)=(-1)^n \sum_{i=0}^{q}%
\frac{P_{i}%
^{(\alpha+n,\alpha+n)}(0)}{\Vert P_{i}^{(\alpha+n,\alpha+n)}\Vert^{2}_{\alpha+n,\alpha+n}}
\frac{d^n}{d t^n} \left[ P_{i}^{(\alpha+n,\alpha+n)}(t) w_{\alpha+n,\alpha+n}(t)\right].
\end{equation}
Since $P_{i}^{(\alpha+n, \alpha+n)}(0)=0$ for any odd integer $i$, (\ref{Q3}) becomes
\begin{equation}\label{Q4}
    Q_{\alpha,n,q}(t)=(-1)^n \sum_{i=0}^{\frac{q}{2}}%
\frac{P_{2i}%
^{(\alpha+n,\alpha+n)}(0)}{\Vert P_{2i}^{(\alpha+n,\alpha+n)}\Vert^{2}_{\alpha+n,\alpha+n}}
\frac{d^n}{d t^n} \left[ P_{2i}^{(\alpha+n,\alpha+n)}(t) w_{\alpha+n,\alpha+n}(t)\right].
\end{equation}
Since $P_{2i}^{(\alpha+n,\alpha+n)}(-t) = P_{2i}^{(\alpha+n,\alpha+n)}(t) $ (see \cite{G.
Szego}) and $w_{\alpha+n,\alpha+n}(-t)=w_{\alpha+n,\alpha+n}(t)$, we have
\begin{equation*}
\frac{d^n}{d t^n} \left[ P_{2i}^{(\alpha+n,\alpha+n)}(t) w_{\alpha+n,\alpha+n}(t)\right]=(-1)^n \frac{d^n}{d t^n} \left[ P_{2i}^{(\alpha+n,\alpha+n)}(-t) w_{\alpha+n,\alpha+n}(-t)\right].
\end{equation*}
Thus, we have $Q_{\alpha,n,q}(t)=(-1)^n Q_{\alpha,n,q}(-t)$. Then this proof is completed. \hfill$\Box$\\

\noindent\textbf{Proof of Proposition \ref{proposition6}.} Let us
recall the local Taylor formula with the Peano remainder term
\cite{Zorich}. For any given $\varepsilon'>0$, there exists
$\delta>0$ such that

\begin{equation}\label{peano1}
\left|f(x+ht)-f_{n-1}(x+ht)-\frac{f_-^{(n)}(x)}{n!}(ht)^n\right|<
\varepsilon' |ht|^n, \text{ for } \delta<ht<0,
\end{equation}
and
\begin{equation}\label{peano2}
\left|f(x+ht)-f_{n-1}(x+ht)-\frac{f_+^{(n)}(x)}{n!}(ht)^n\right|<
\varepsilon' (ht)^n,\text{ for } 0 <ht<\delta,
\end{equation}
where $f_{n-1}(x+ht)$ is the $(n-1)${th} order truncated Taylor
series expansion of $f(x+ht)$. Let us consider the function
$g(x)=x^n$  the $n${th} order derivative of which is equal to
$(n!)$. Thus, by using (\ref{d_h_N}) we have
\begin{equation*}
\forall x\in I_{h},\ D_{h,\alpha,\alpha,q}^{(n)}g(x)=(n!).
\end{equation*}
Thus, by applying Corollary \ref{proposition2} with $\alpha=\beta$,
we get
\begin{equation*}
\forall x\in I_{h},\
D_{h,\alpha,\alpha,q}^{(n)}g(x)=\frac{1}{h^{n}}\int
_{-1}^{1}Q_{\alpha,n,q}(t)\,g(x+ht)\,dt=(n!).
\end{equation*}
In particular, by taking  $x=0$ we get $ \frac{1}{h^{n}}\int
_{-1}^{1}Q_{\alpha,n,q}(t)\,(ht)^n\,dt=(n!). $ According to Lemma
\ref{lemma3}, $t^n Q_{\alpha,n,q}(t)$ with $t \in [-1,1]$ is an odd
function. Hence, we have $ \frac{1}{h^{n}}\int
_{-1}^{0}Q_{\alpha,n,q}(t)\,(ht)^n\,dt=\frac{1}{h^{n}}\int
_{0}^{1}Q_{\alpha,n,q}(t)\,(ht)^n\,dt. $ Thus, we get
\begin{equation} \label{52}
\frac{1}{h^{n}}\int
_{-1}^{0}Q_{\alpha,n,q}(t)\frac{f_-^{(n)}(x)}{n!}(ht)^n\,dt=\frac{1}{2}
f_-^{(n)}(x),
\end{equation}
and
\begin{equation} \label{53}
\frac{1}{h^{n}}\int
_{0}^{1}Q_{\alpha,n,q}(t)\frac{f_+^{(n)}(x)}{n!}(ht)^n\,dt=\frac{1}{2}
f_+^{(n)}(x).
\end{equation}
By using (\ref{d_h_N}) and Corollary \ref{proposition2} with
$\alpha=\beta$ we get
\begin{equation} \label{54}
\forall x\in I_{h},\
D_{h,\alpha,\alpha,q}^{(n)}f_{n-1}(x)=\frac{1}{h^{n}}\int
_{-1}^{1}Q_{\alpha,n,q}(t)\,f_{n-1}(x+ht)\,dt=0.
\end{equation}
Hence, by using (\ref{52}), (\ref{53}) and (\ref{54})   we obtain
\begin{equation}
\begin{split}\label{inegalite}
&\left|
 D_{h,\alpha,\alpha,{q}}^{(n)}f(x)-\frac{1}{2}\left(f^{(n)}_+(x)+f^{(n)}_-(x)\right)\right|\\
 \leq & \frac{1}{h^{n}}\int
_{-1}^{0} \left|Q_{\alpha,n,q}(t)\left(
f(x+ht)-f_{n-1}(x+ht)-\frac{f_-^{(n)}(x)}{n!}(ht)^n\right) \right|dt
\\+& \frac{1}{h^{n}}\int _{0}^{1} \left|Q_{\alpha,n,q}(t)\left(
f(x+ht)-f_{n-1}(x+ht)-\frac{f_+^{(n)}(x)}{n!}(ht)^n\right)
\right|dt.
\end{split}
\end{equation}
By using (\ref{Q2}), we get
\begin{equation}
\int _{0}^{1}\left|Q_{\alpha,n,q}(t)\right|dt \leq
\displaystyle\sum_{i=0}^{\frac{q}{2}}{\left|P_{2i}^{(\alpha+n,\alpha
+n)}(0)\right|}\sum_{j=0}^{2i}\binom{2i}{j}\frac{4i+2\alpha
+2n+1}{2i+2\alpha+2n+1}\int
_{0}^{1}\left|\rho_{n,\alpha_{2i,j},\beta_{j}}(t)\right|dt.
\end{equation}
Then, according to (\ref{d_h}) and (\ref{jacobi_w}) we can obtain
that $\int
_{0}^{1}\left|\rho_{n,\alpha_{2i,j},\beta_{j}}(t)\right|dt <
\infty$. Hence, $$\int _{0}^{1}\left|Q_{\alpha,n,q}(t)\,t^n\right|dt
\leq \int _{0}^{1}\left|Q_{\alpha,n,q}(t)\right|dt < \infty.$$
Consequently, for any $\varepsilon>0$, by using (\ref{inegalite}),
(\ref{peano1}) and (\ref{peano2}) with $\varepsilon=2 \varepsilon'
\int _{0}^{1}\left|Q_{\alpha,n,q}(t)\,t^n\right|dt$, there exists $\delta$ such that
$0<h<\delta$ and
\begin{equation*}%\label{}
\left|
 D_{h,\alpha,\alpha,{q}}^{(n)}f(x)-\frac{1}{2}\left(f^{(n)}_+(x)+f^{(n)}_-(x)\right)\right|<\varepsilon.
\end{equation*}
Then, this proof can be completed. \hfill $\Box$
%%%%%%%%%%%%%%%%%%%%%%%%%%%%%%%%%%%%%%%%%%%%%%%%%%%%%%%%%%%%%%%%%%%%%%%%%%%%%%%%%%%%%%%%%%%%%%%%%%%%%%%%%%%%%%

\begin{figure}[h!]
 \centering
 \subfigure[$n=1$, $q=0,2,\cdots,8$ and
$\alpha=0,1,\cdots,10$.]
 {\includegraphics[scale=0.55]{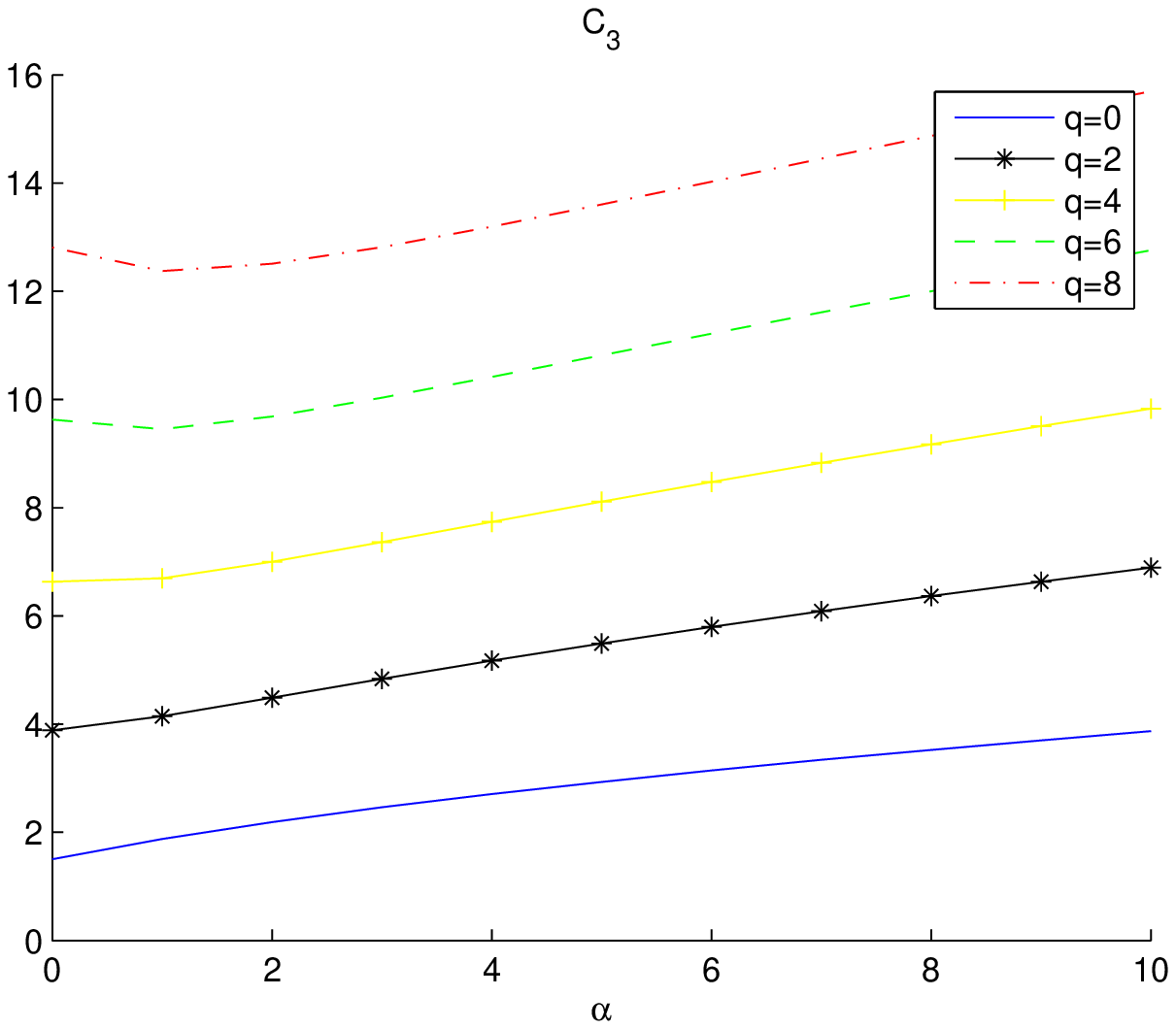}}
 \subfigure[ $n=2$, $q=0,2,\cdots,8$ and
$\alpha=0,1,\cdots,10$.]
 {\includegraphics[scale=0.55]{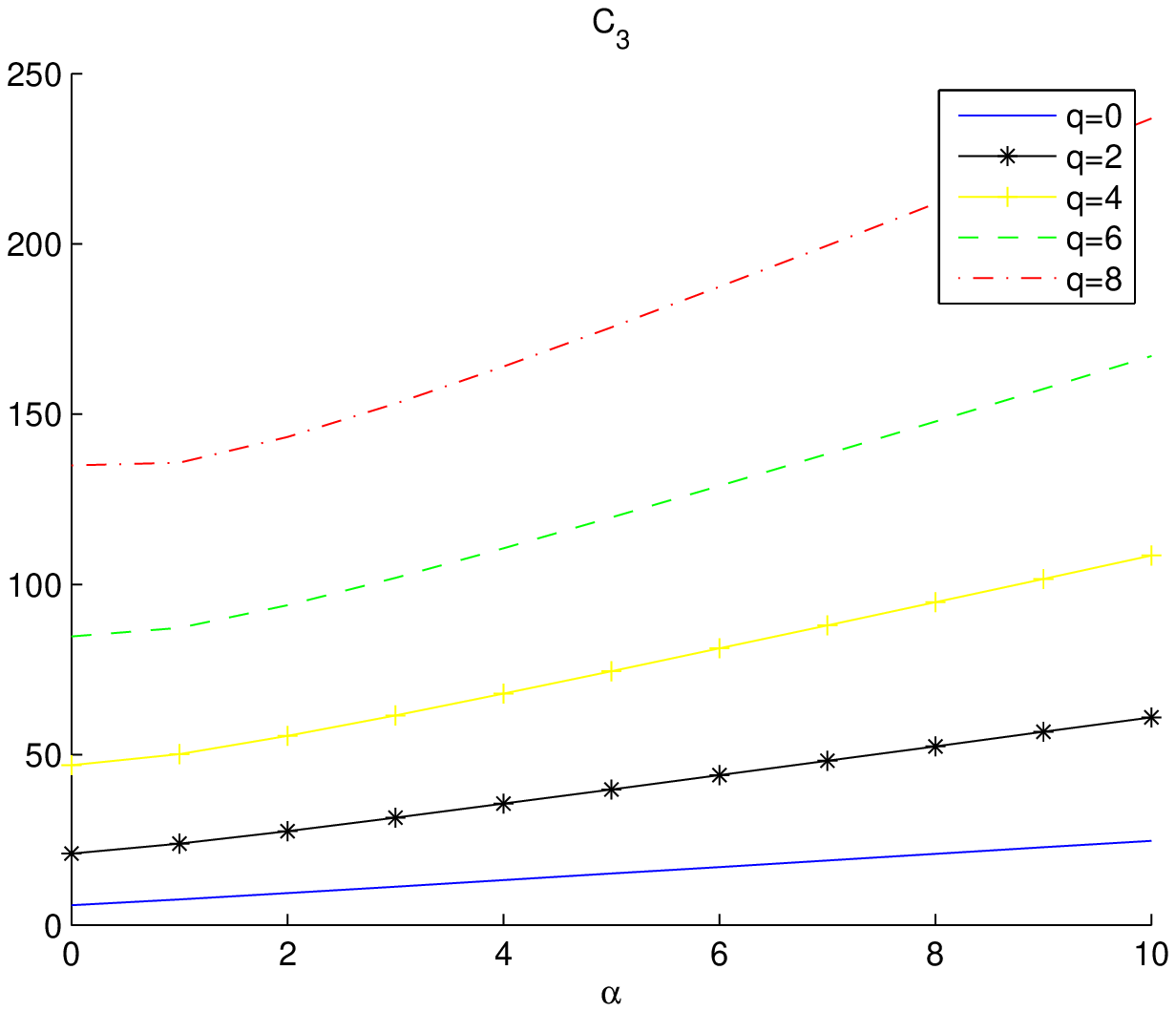}}%\label{figureq2}
\caption{ Values of $C_{3}$.}\label{figureq}
\end{figure}

\begin{figure}[h!]
 \centering
 \subfigure[$q=4$, $n=1,\cdots,4$ and
$\alpha=0,1,\cdots,10$.]
 {\includegraphics[scale=0.55]{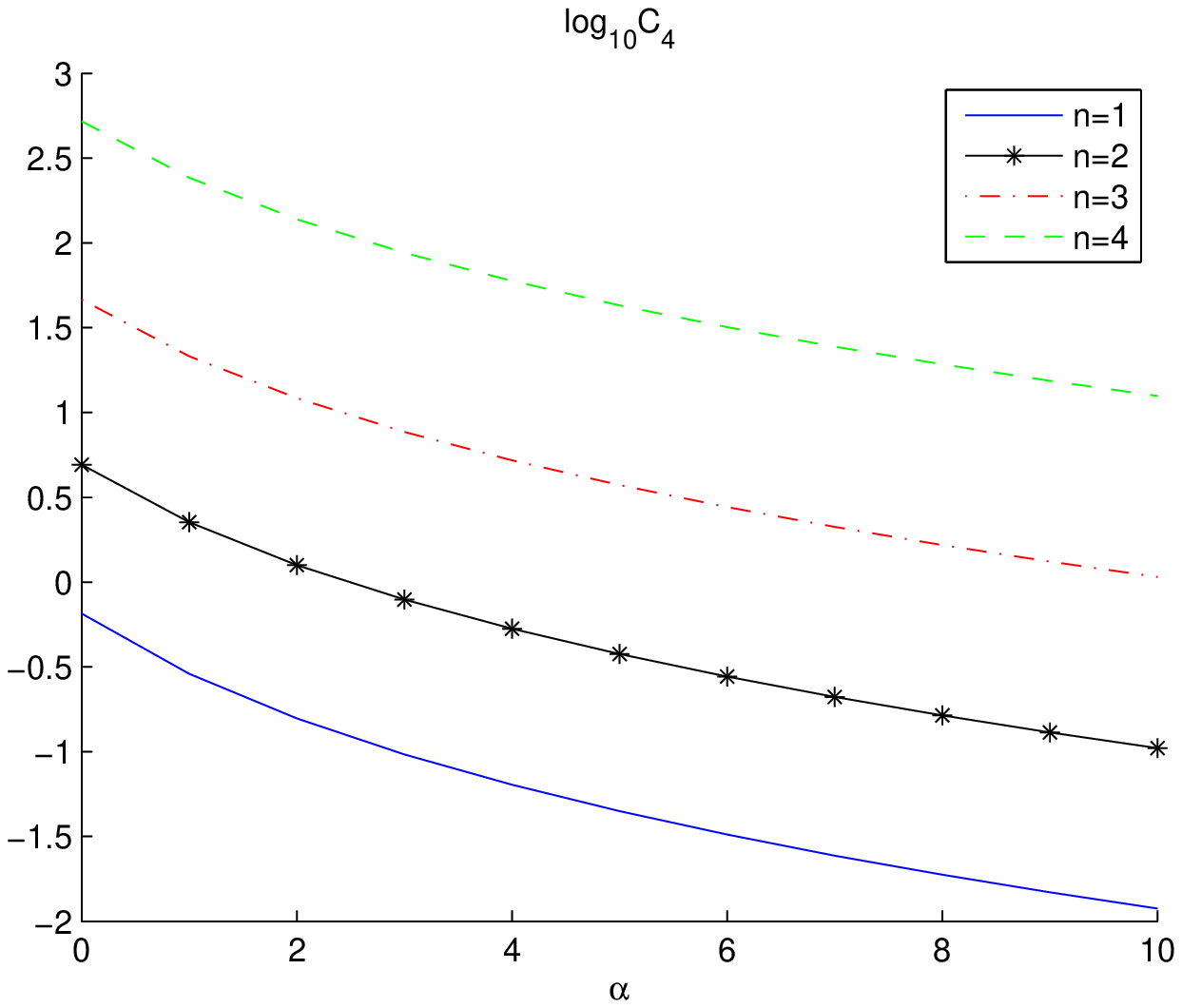}}
 \subfigure[  $q=4$, $n=1,\cdots,4$ and
$\alpha=0,1,\cdots,10$.]
 {\includegraphics[scale=0.55]{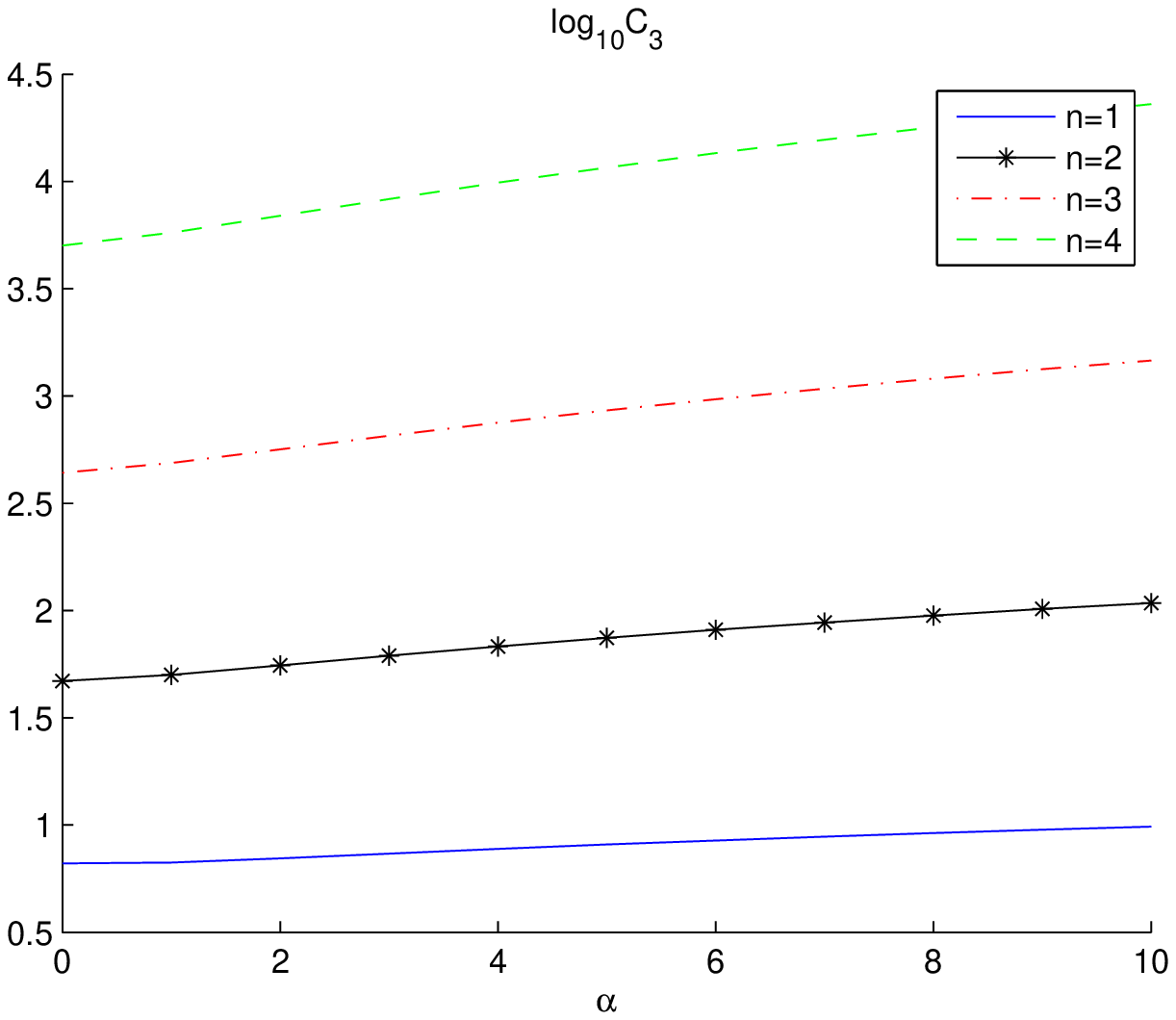}}
\caption{ Values of $\log_{10}C_{4}$ and
$\log_{10}{C}_{3}$.}\label{figurealpha}
\end{figure}

\section{Numerical tests} \label{section3}

In order to demonstrate the efficiency  and the  stability of the
previously proposed estimators, we present some numerical results in
this section. First of all, we analyze the choice of parameters for
these estimators.

\subsection{Analysis for parameters' choice for the bias term error and the noise error}
As  is shown previously, the proposed estimators contain two
sources of errors: the bias term error which is produced by the
truncation of the Jacobi orthogonal series expansion and the noise
error contribution. The error bounds for these errors are given in
Corollary \ref{corollary5}. We are going to use the knowledge of the
parameters' influence to these error bounds. This will help us to
obtain a tendency on  the influence of these parameters on the
estimation errors.

According to Corollary \ref{corollary4}, we set $\alpha=\beta$ and
choose the truncation order $q$ to be an even integer. On the one hand,
it is clear that we should set $q$ as large as possible so as to
improve the convergence rate and reduce the bias term error. On the
other hand, the noise error contribution is bounded in Corollary
\ref{corollary5} by $B_{noise}=C_{3}\frac{\delta}{h^{n}}$ where
$C_{3}=\int_{-1}^{1} | Q_{\alpha,n,q}(t)| \, dt$. We can see in
Figure \ref{figureq} the different values  of $C_{3}$ where $n=1,2$,
$q=0,2,\cdots,8$ and $\alpha=0,1,\cdots,10$. It is clear that with
the same values for $n$ and $\alpha$, $C_{3}$  increases with
respect to $q$. Furthermore, it is easy to verify that  $C_{3}$
increases with respect to $q$, independently of $n$ and $\alpha$.
Hence, in order to reduce the bias term error and to avoid  a large
noise error, we set $q=4$ in our estimators. With this value,
according to Corollary \ref{corollary5} the convergence rate is
$O(\delta^{\frac{6}{n+6}})$.

The bias term error is bounded by $B_{bias}=\hat{C}_{2}h^{q+2}$ in
Corollary
\ref{corollary5} where $\hat{C}_{2}=\frac{{M}_{n+2+q}}{(n+q+2)!}\int_{-1}^{1}|t^{n+q+2}%
Q_{\alpha,n,q}(t)|\,dt$. Let us introduce $C_4=\int_{-1}^{1}|t^{n+q+2}%
Q_{\alpha,n,q}(t)|\,dt$. We can see in Figure \ref{figurealpha} the
different values of $\log_{10}C_{4}$ and $\log_{10}C_3$ when
$n=1,\cdots,4$, $q=4$ and $\alpha=0,1,\cdots,10$. It is clear that
$C_4$ decreases  with respect to $\alpha$ while ${C}_3$ increases
with respect to $\alpha$. Thus, in order to reduce the bias term
error, we should set $\alpha$ as large as possible. However, a large
value of $\alpha$ may produce a large noise error contribution.
Here, we choose $\alpha=5$.

Until here, we have chosen $q=4$ and $\alpha=5$. The noise error
decreases with respect to $h$ and the bias term error  increases
with respect to $h$.  In the next subsection we are going to choose
an appropriate value for $h$ by using the knowledge of function $f$
and by taking into account the numerical integration method error.

%\clearpage

\begin{figure}[h!]
 \centering
 \subfigure[$n=1$, $\alpha=5$, $q=4$ and $h=591T_s$.]
 {\includegraphics[scale=0.55]{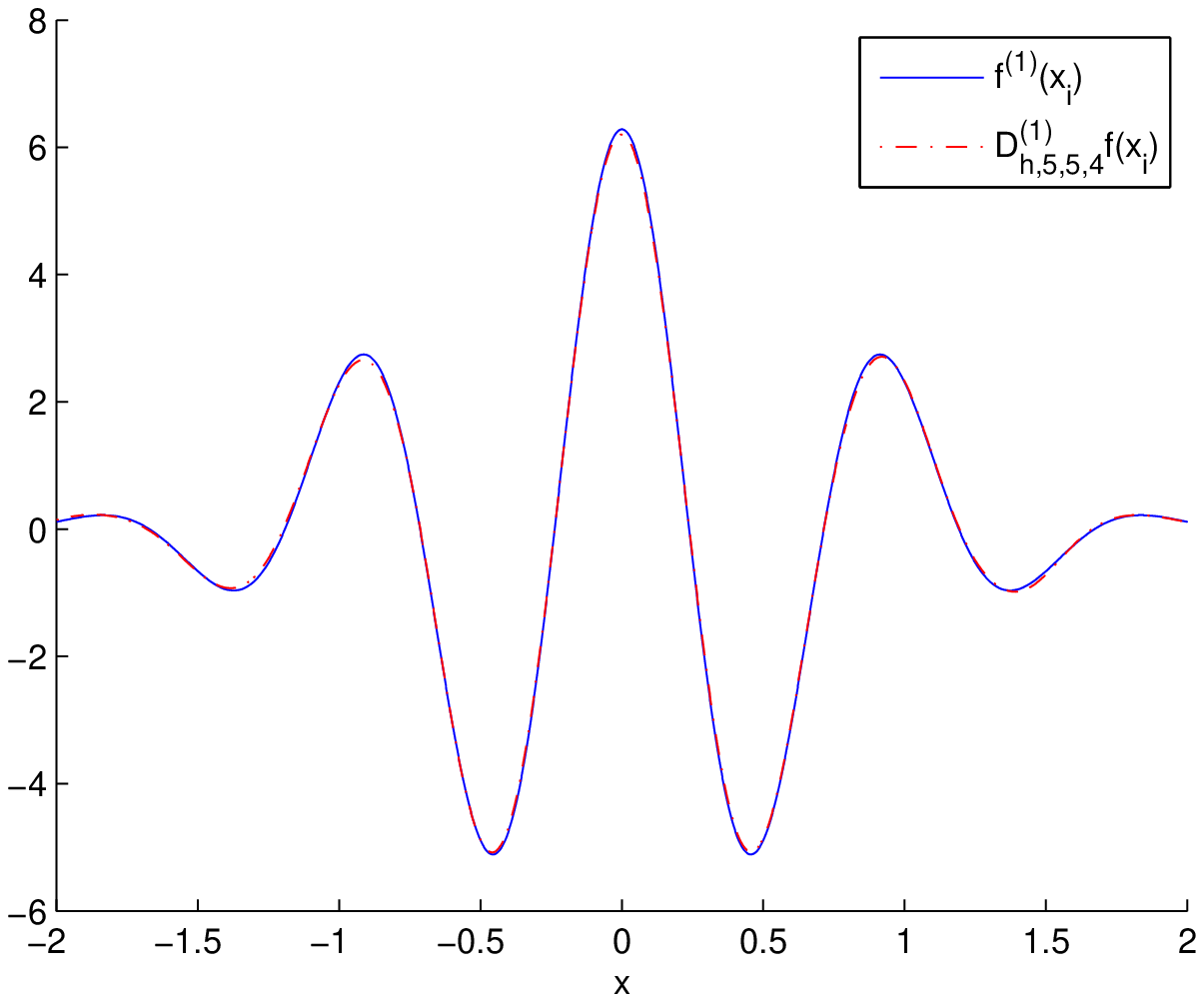}}
 \subfigure[ $n=2$, $\alpha=5$, $q=4$ and $h=698T_s$.]
 {\includegraphics[scale=0.55]{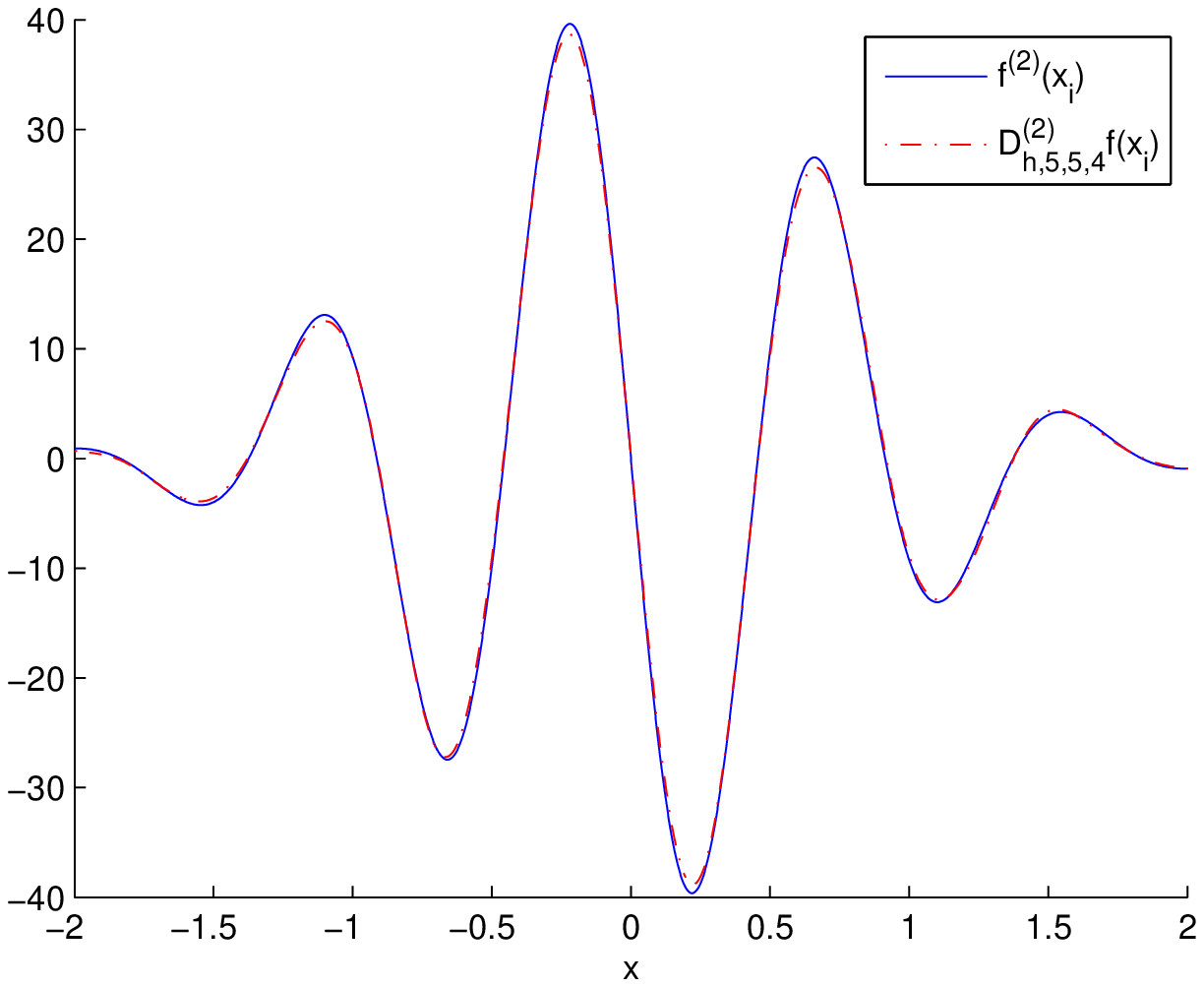}}
 \subfigure[ $n=3$, $\alpha=5$, $q=4$ and $h=777T_s$.]
 {\includegraphics[scale=0.55]{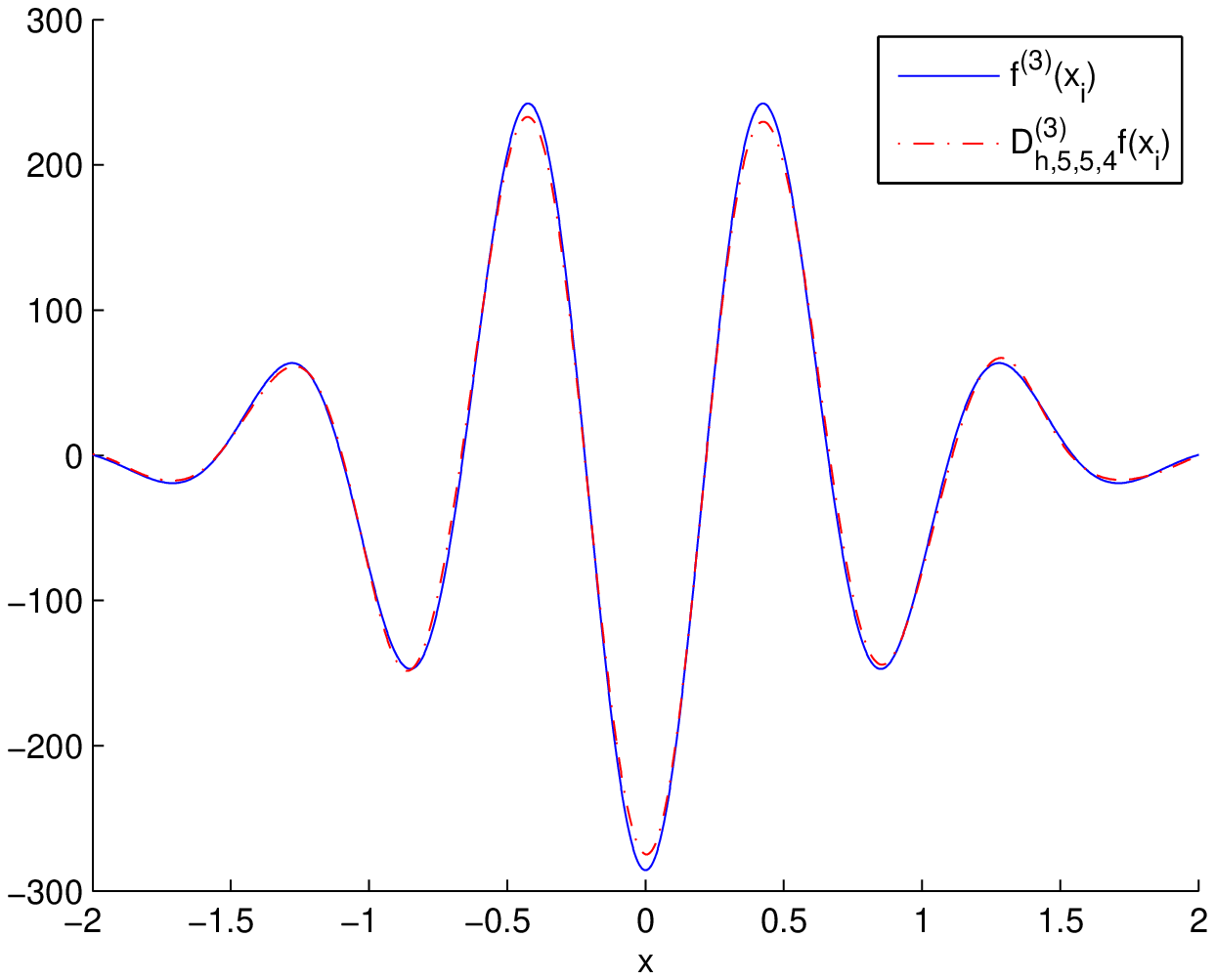}}
  \subfigure[ $n=4$, $\alpha=5$, $q=4$ and $h=850T_s$.]
 {\includegraphics[scale=0.55]{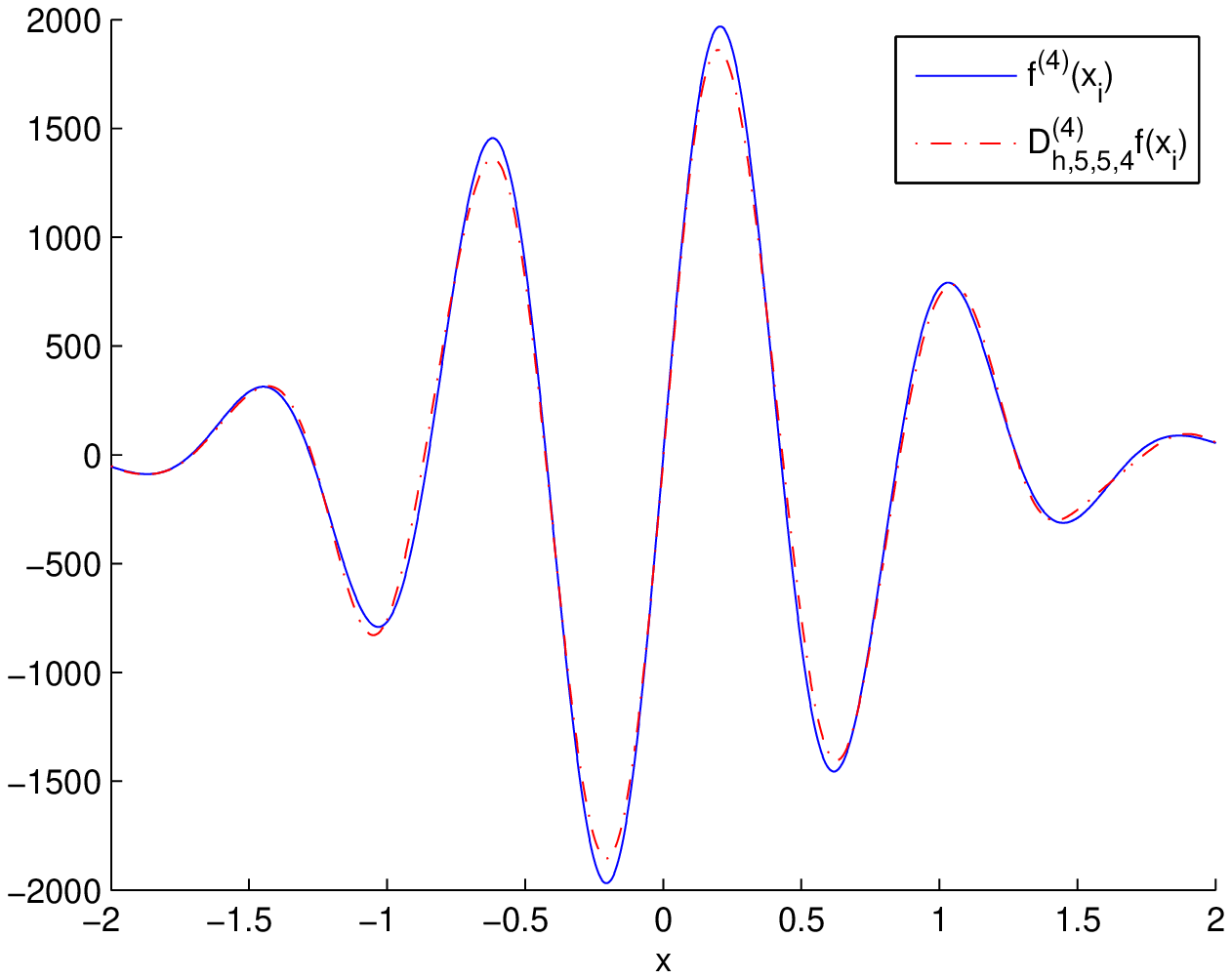}}
\caption{The exact values of  $f^{(n)}_1(x_i)$ and the estimated
values $D_{h,\alpha,\alpha,q}^{(n)}f_1(x_i)$ for $\delta=0.15$.}
\label{fig2}
\end{figure}

\begin{figure}[h!]
 \centering
 \subfigure[$n=1$, $\alpha=5$, $q=4$ and $h=442T_s$.]
 {\includegraphics[scale=0.55]{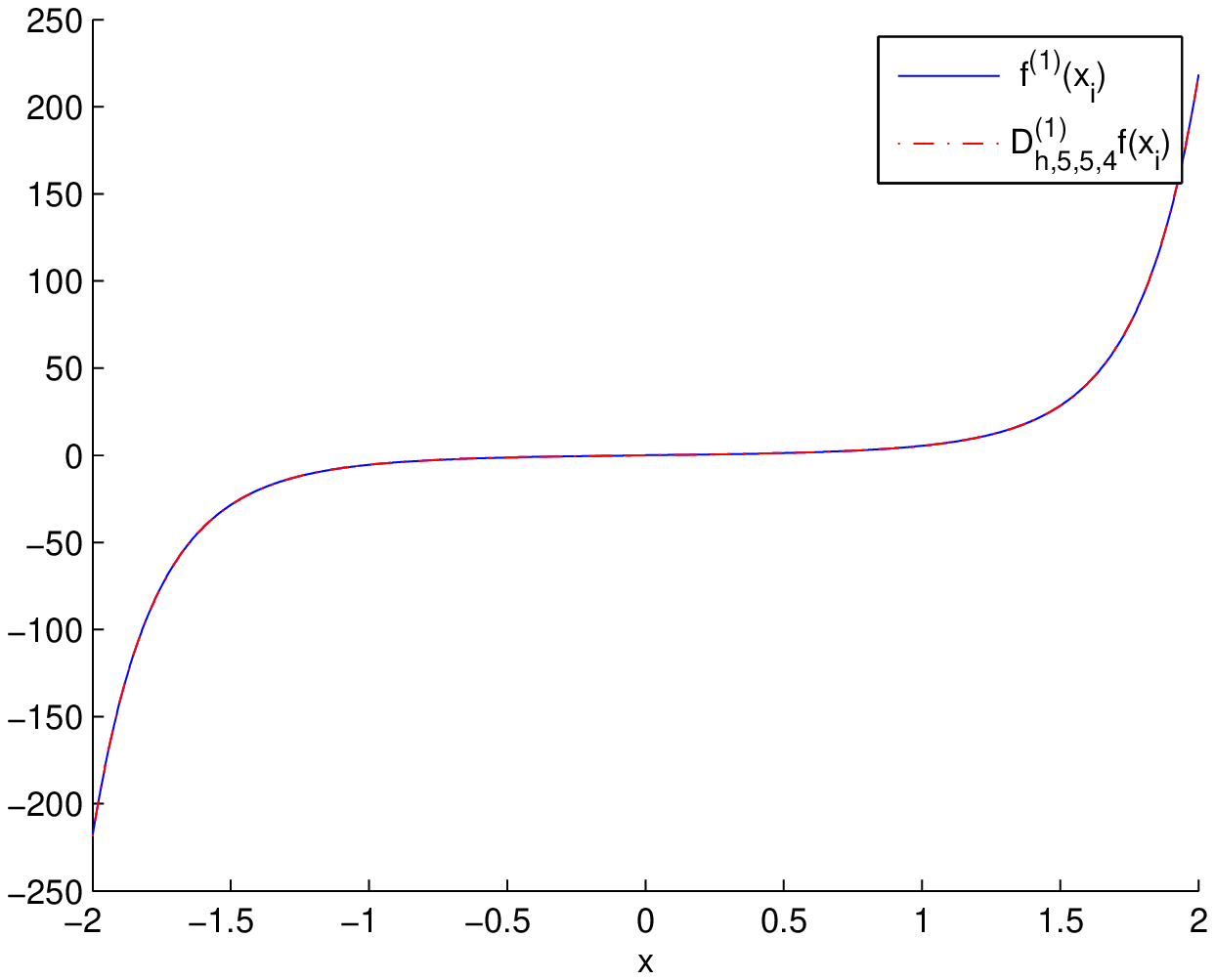}}
 \subfigure[ $n=2$, $\alpha=5$, $q=4$ and $h=549T_s$.]
 {\includegraphics[scale=0.55]{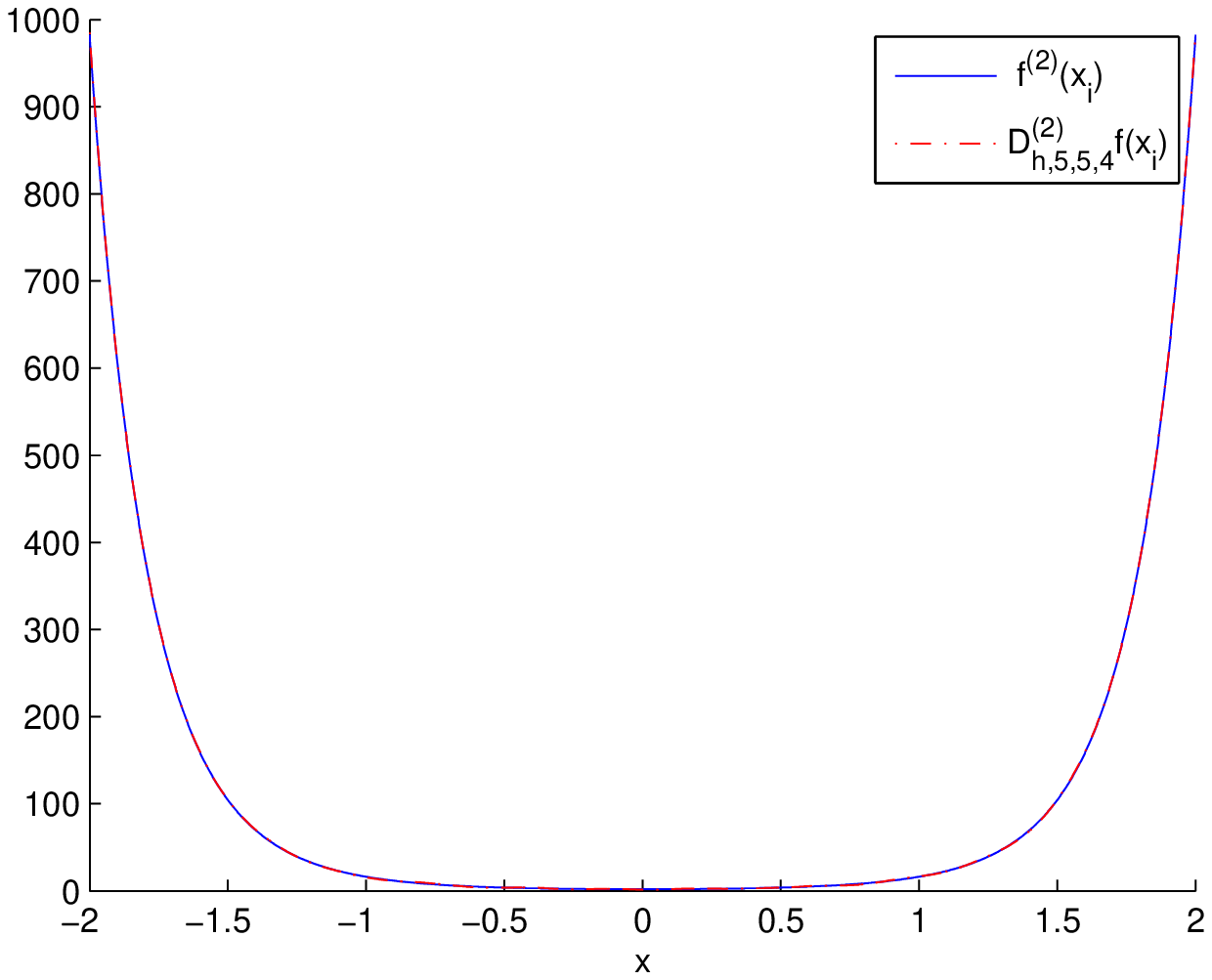}}
 \subfigure[ $n=3$, $\alpha=5$, $q=4$ and $h=643T_s$.]
 {\includegraphics[scale=0.55]{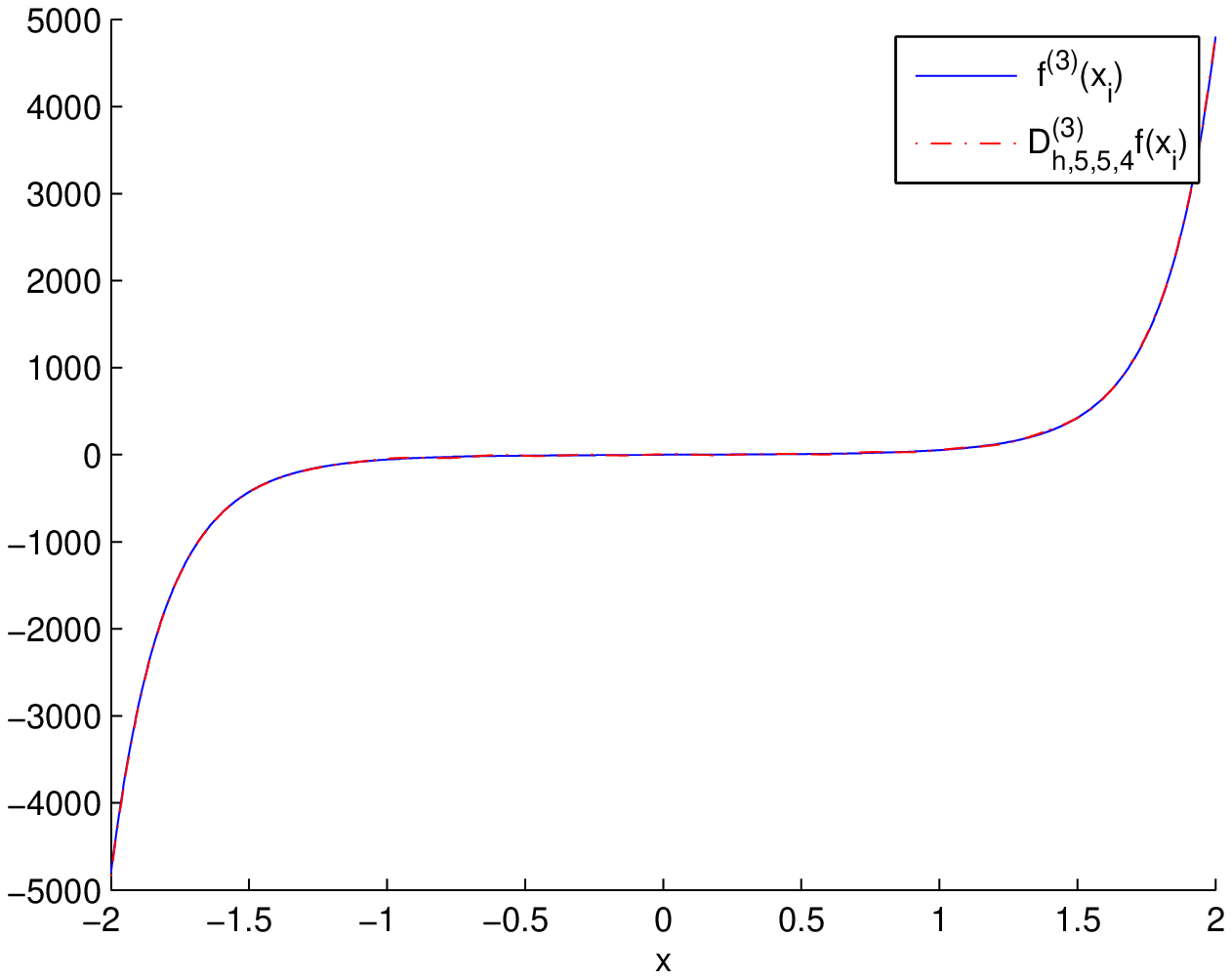}}
  \subfigure[ $n=4$, $\alpha=5$, $q=4$ and $h=733T_s$.]
 {\includegraphics[scale=0.55]{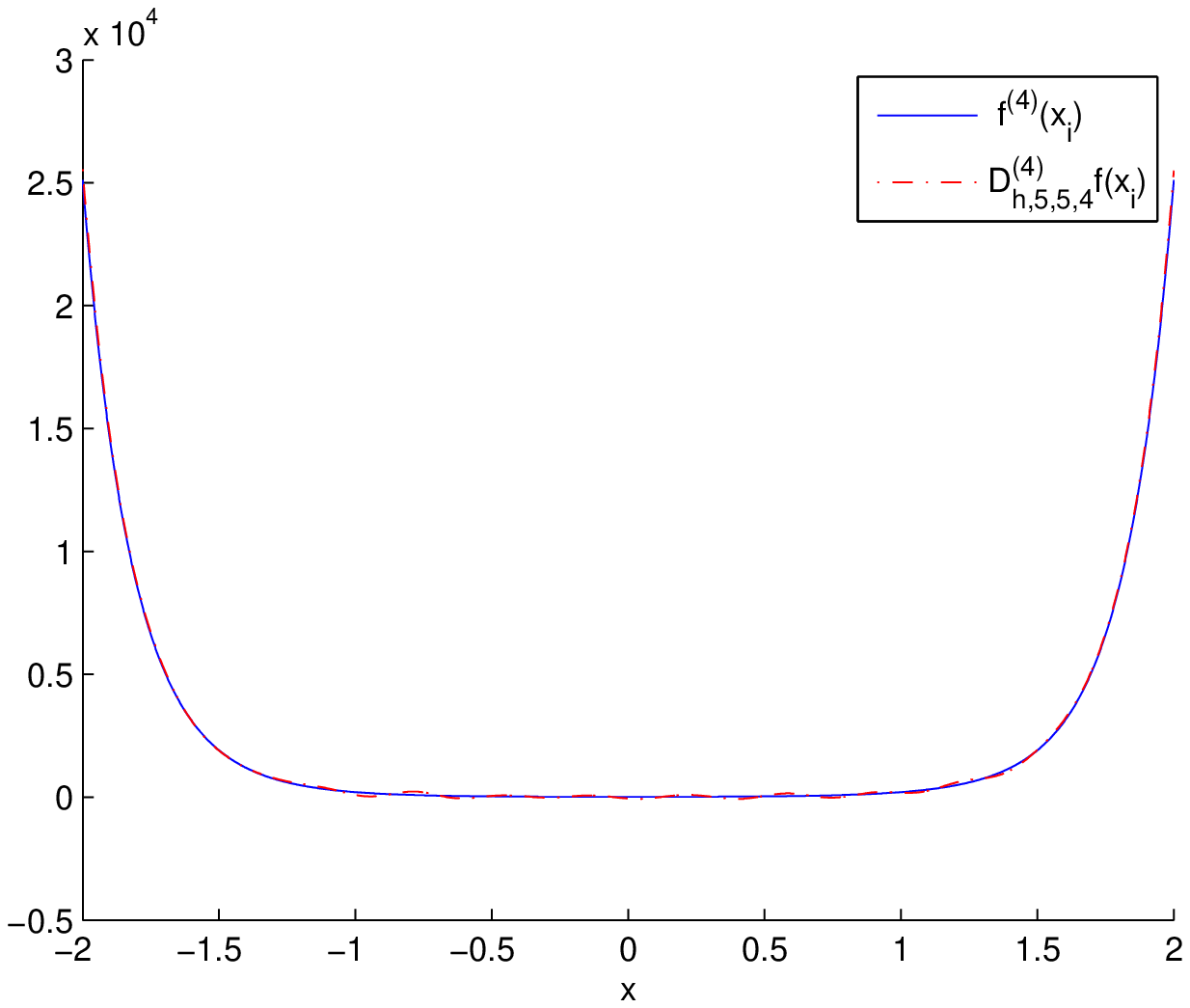}}
\caption{The exact values of  $f^{(n)}_2(x_i)$ and the estimated
values $D_{h,\alpha,\alpha,q}^{(n)}f_2(x_i)$ for $\delta=0.15$.}
\label{fig1}
\end{figure}

\begin{figure}[h!]
 \centering
 \subfigure[$n=1$, $\alpha=5$, $q=4$ and $h=442T_s$.]
 {\includegraphics[scale=0.55]{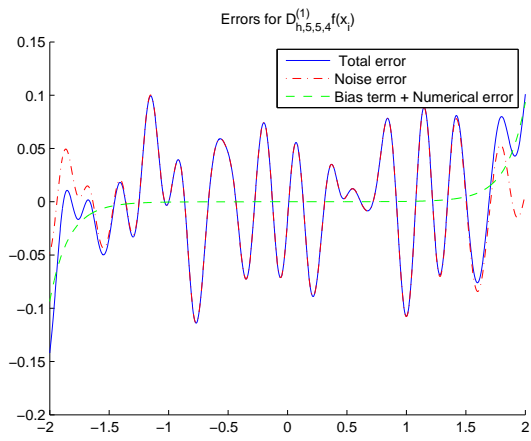}}
 \subfigure[ $n=2$, $\alpha=5$, $q=4$ and $h=549T_s$.]
 {\includegraphics[scale=0.55]{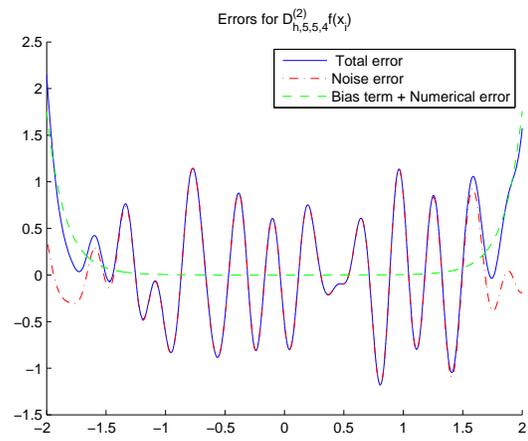}}
 \subfigure[ $n=3$, $\alpha=5$, $q=4$ and $h=643T_s$.]
 {\includegraphics[scale=0.55]{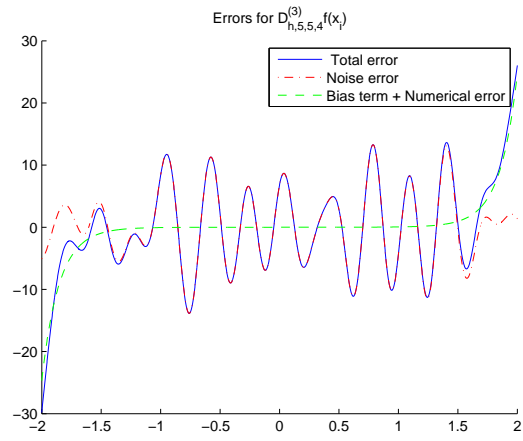}}
  \subfigure[ $n=4$, $\alpha=5$, $q=4$ and $h=733T_s$.]
 {\includegraphics[scale=0.55]{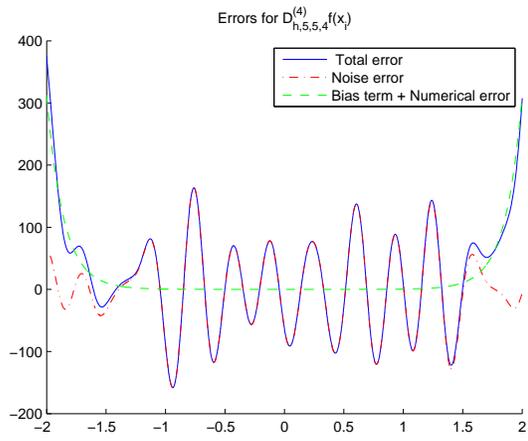}}
\caption{The estimation errors for the estimated values
$D_{h,\alpha,\alpha,q}^{(n)}f_2(x_i)$ for $\delta=0.15$.}
\label{figerror}
\end{figure}

\begin{figure}[h!]
 \centering
\subfigure[The estimation errors for $D_{h_i,5,5,4}^{(n)}f_2(x_i)$
with varying values of $h_i$ for $\delta=0.15$.]
 {\includegraphics[scale=0.55]{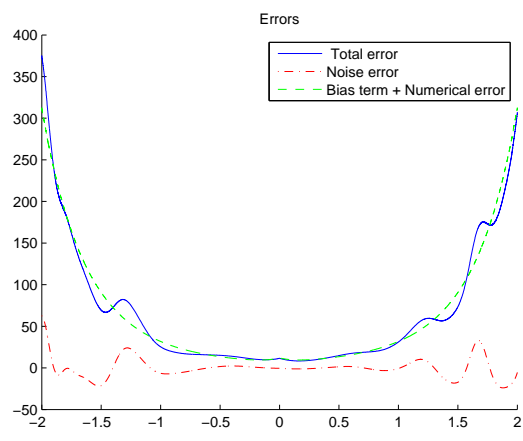}}
\subfigure[Values of $m_i$]
 {\includegraphics[scale=0.55]{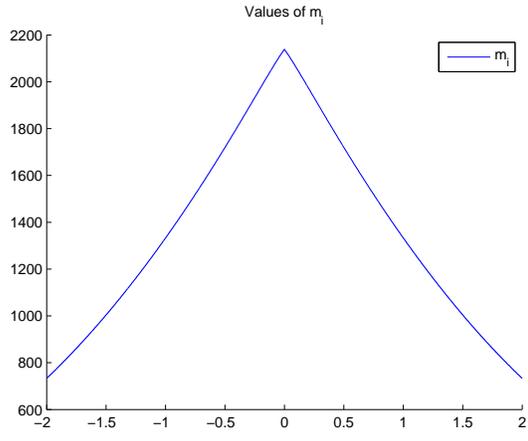}}
\caption{Errors for improved estimations with varying values of $h_i$.}\label{fig11}
\end{figure}

\begin{figure}[h!]
 \centering
 \subfigure[Noise error with its error bounds]
 {\includegraphics[scale=0.55]{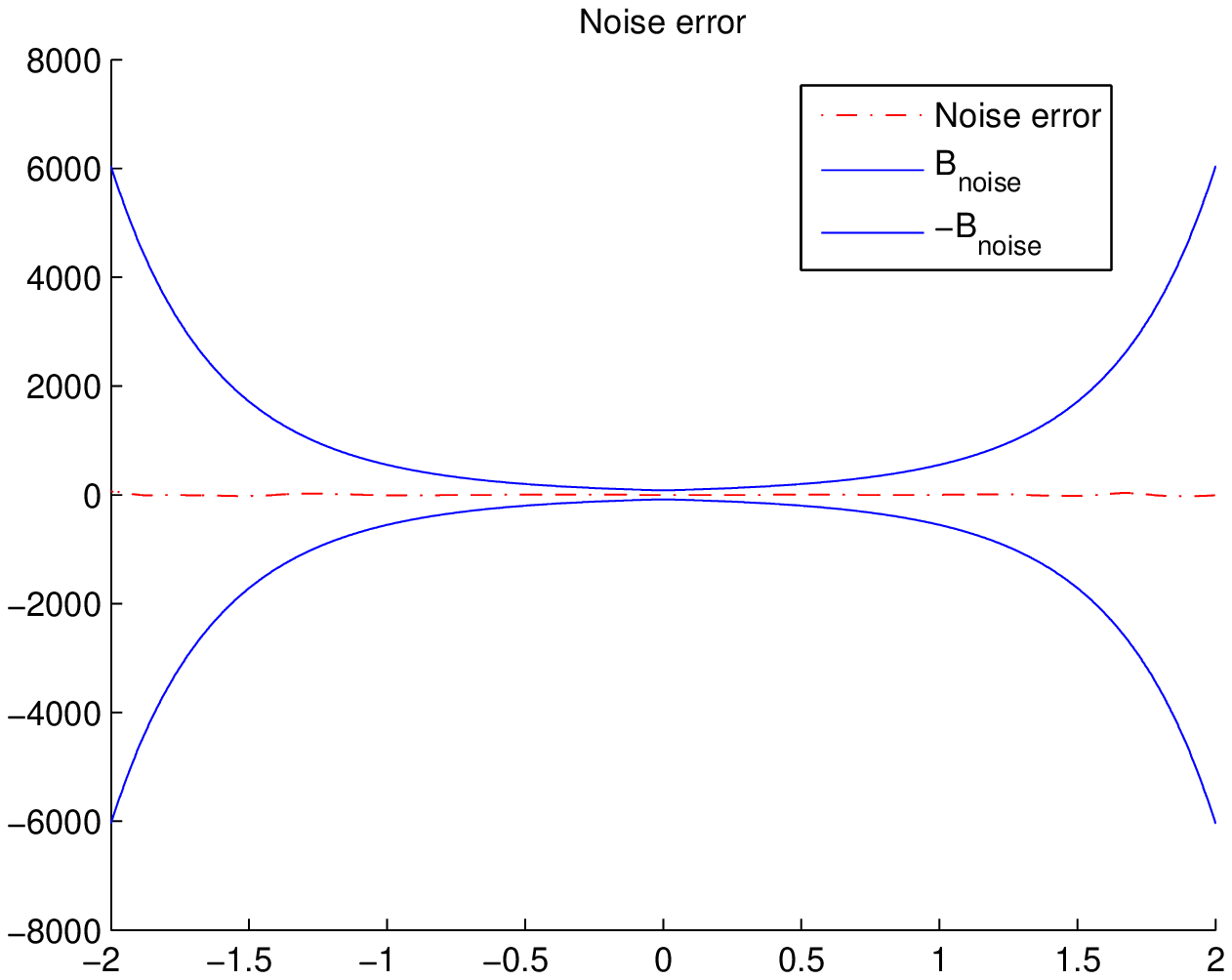}}
 \subfigure[Bias term error + numerical integration error with the error bounds]
 {\includegraphics[scale=0.55]{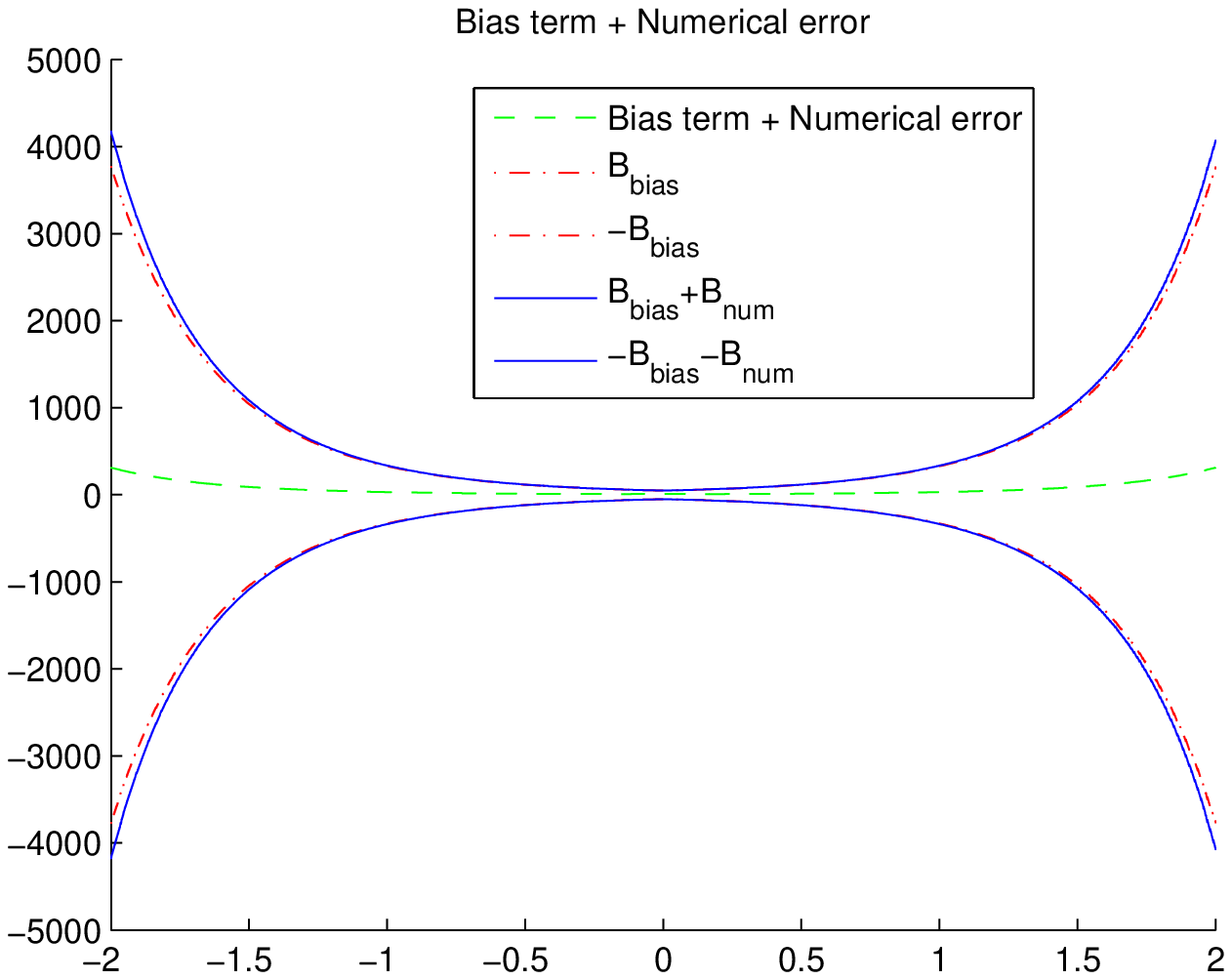}}
\caption{The estimation errors and their corresponding error bounds
for $D_{h_i,5,5,4}^{(n)}f_2(x_i)$ with varying values of $h_i$ for
$\delta=0.15$.} \label{fig22}
\end{figure}

\begin{figure}[h!]
 \centering
 \subfigure[$\delta=0.15$, $n=1$, $\alpha=5$, $q=4$ and $h=1700T_s$.]
 {\includegraphics[scale=0.55]{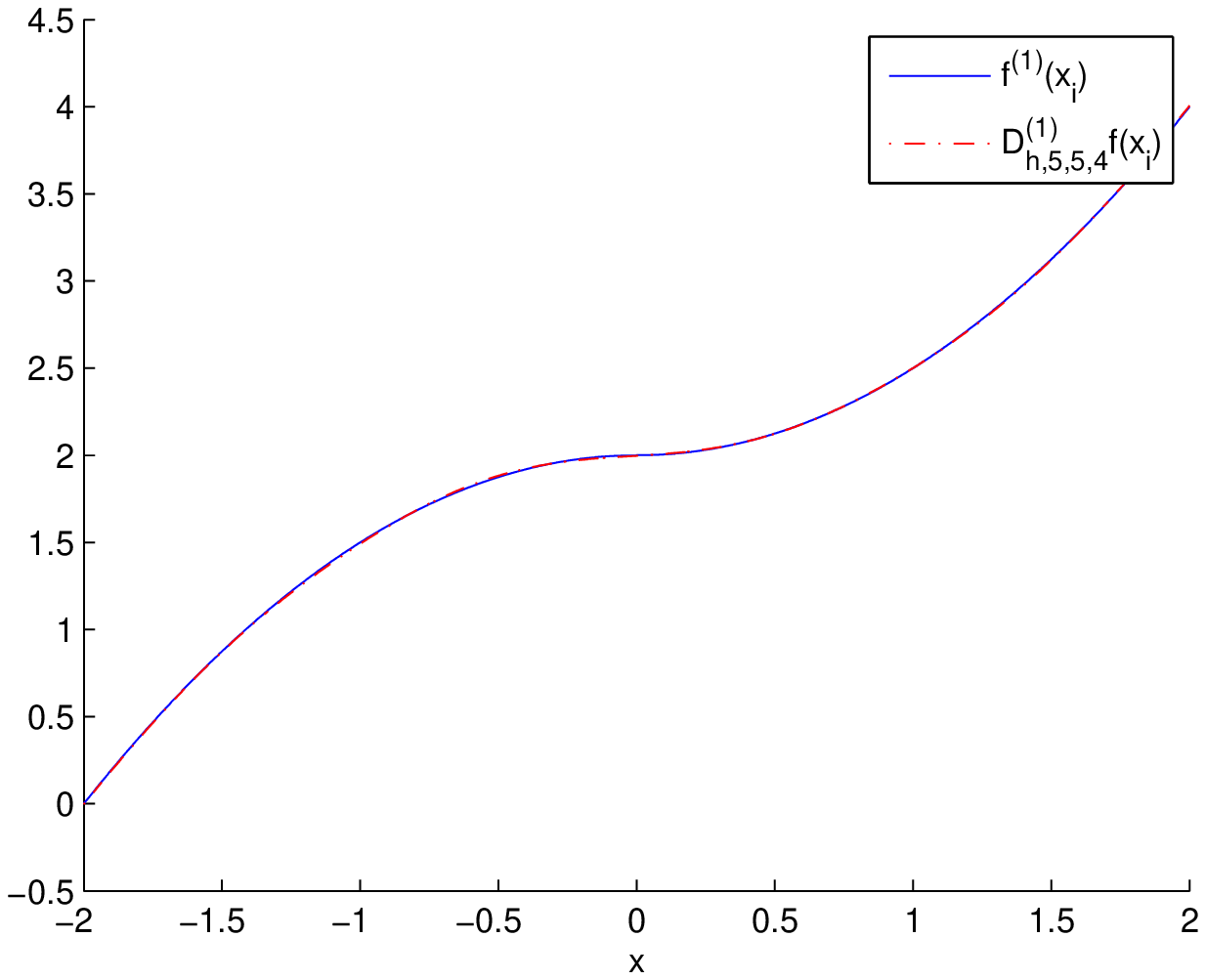}}
  \subfigure[$\delta=0.015$, $n=1$, $\alpha=5$, $q=4$ and $h=1200T_s$.]
 {\includegraphics[scale=0.55]{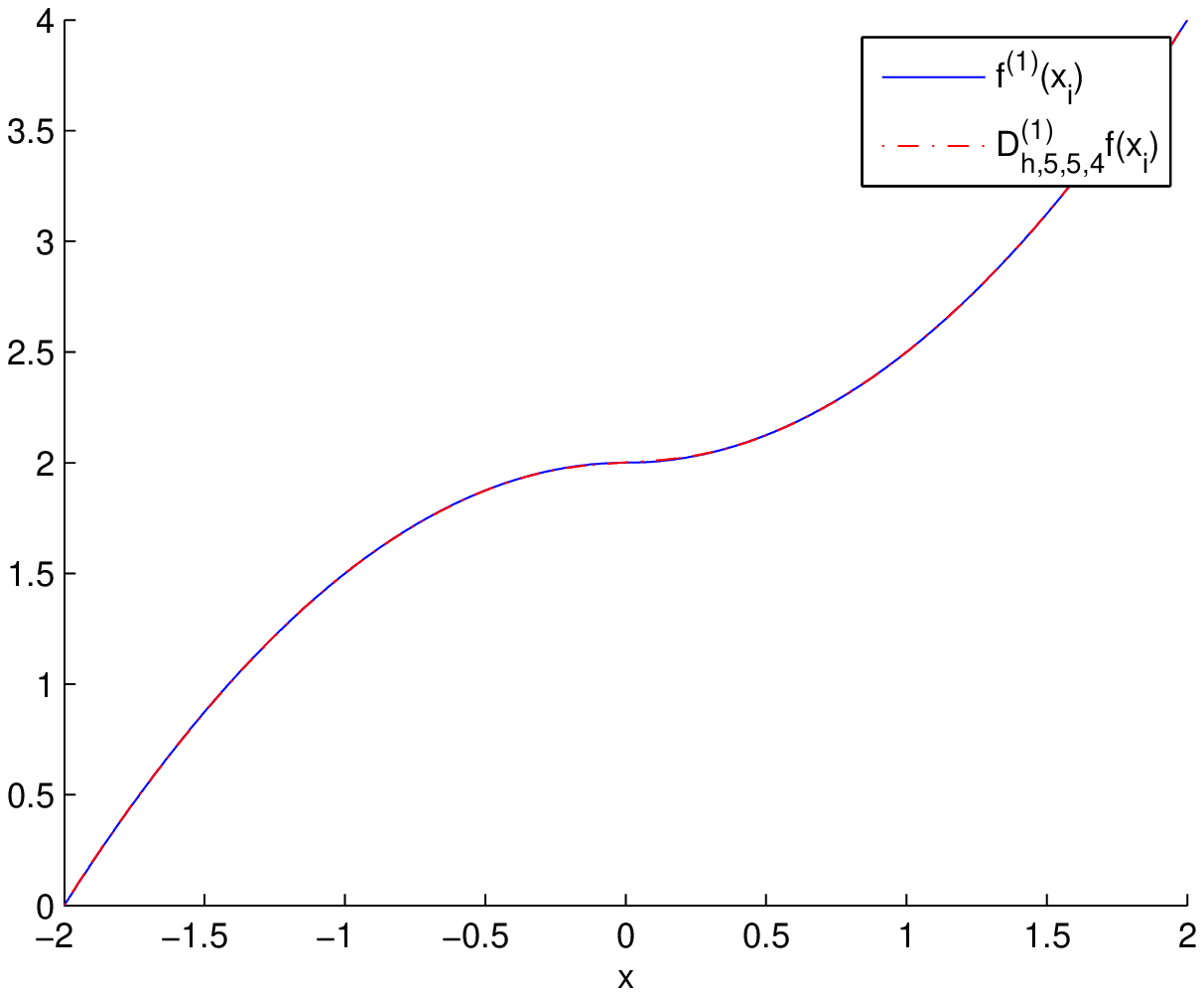}}
 \subfigure[$\delta=0.15$, $n=2$, $\alpha=5$, $q=4$ and $h=1700T_s$.]
 {\includegraphics[scale=0.55]{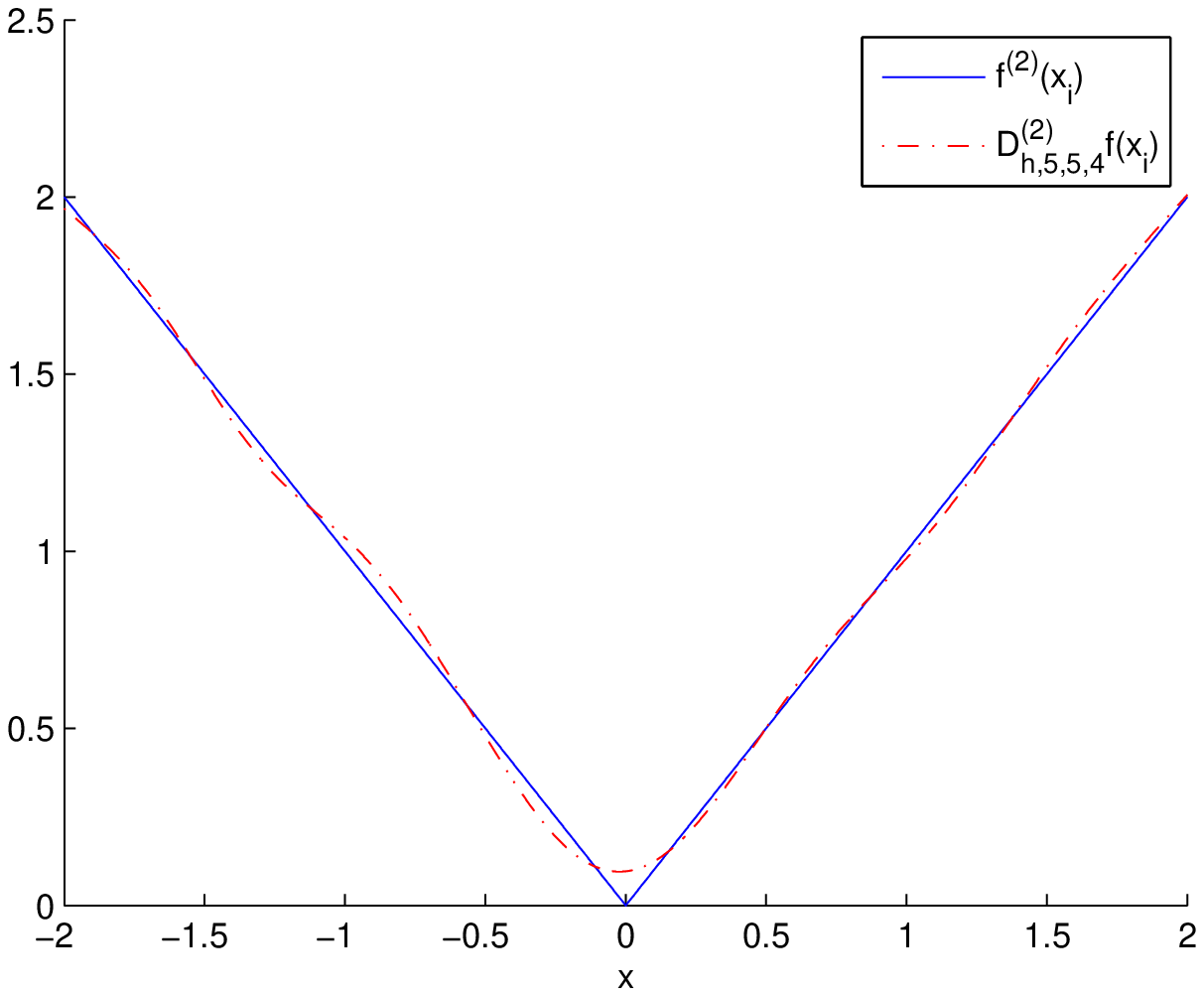}}
  \subfigure[$\delta=0.015$, $n=2$, $\alpha=5$, $q=4$ and $h=1200T_s$.]
 {\includegraphics[scale=0.55]{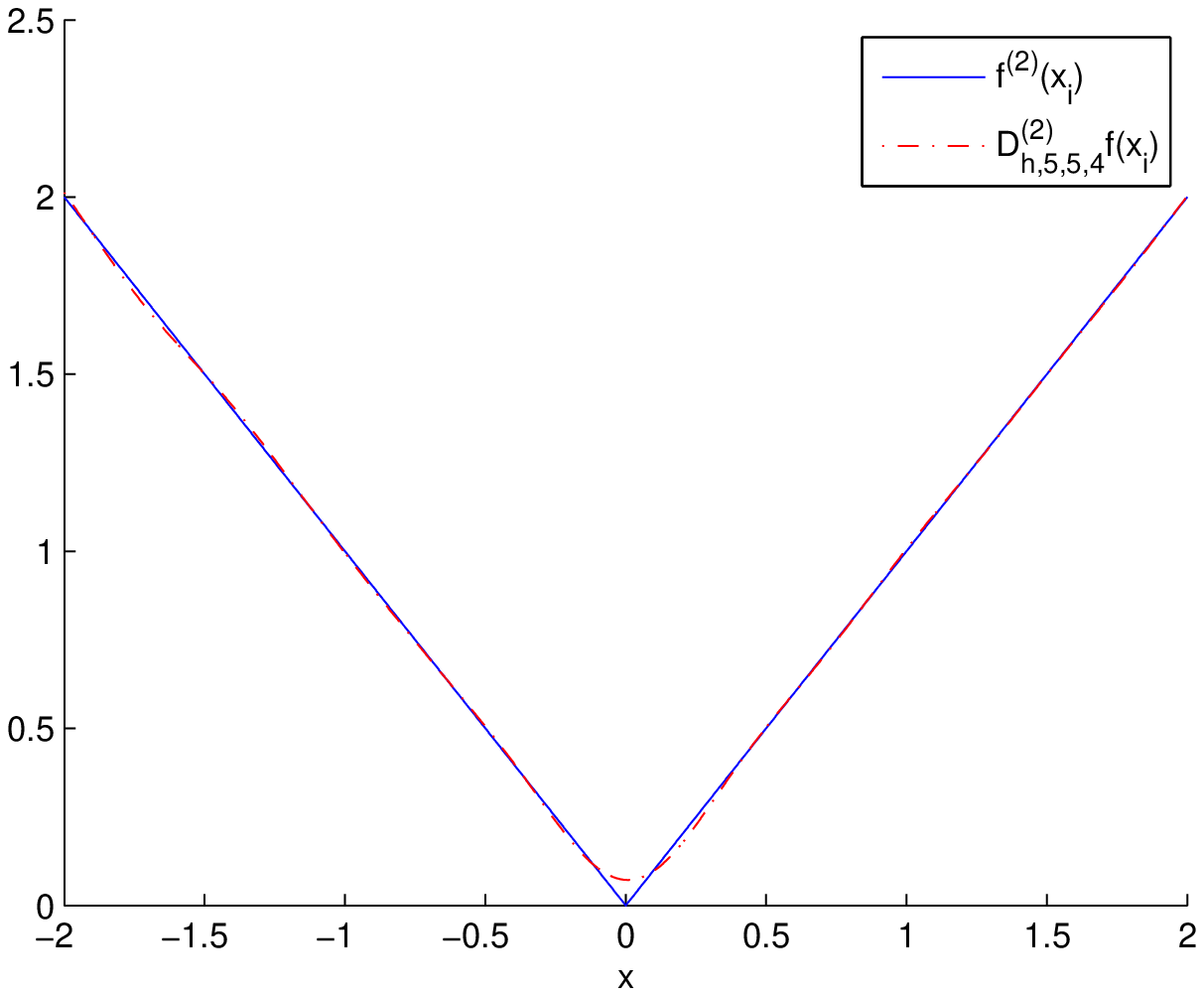}}
  \subfigure[$\delta=0.15$, $n=3$, $\alpha=2$, $q=2$ and $h=1700T_s$.]
 {\includegraphics[scale=0.55]{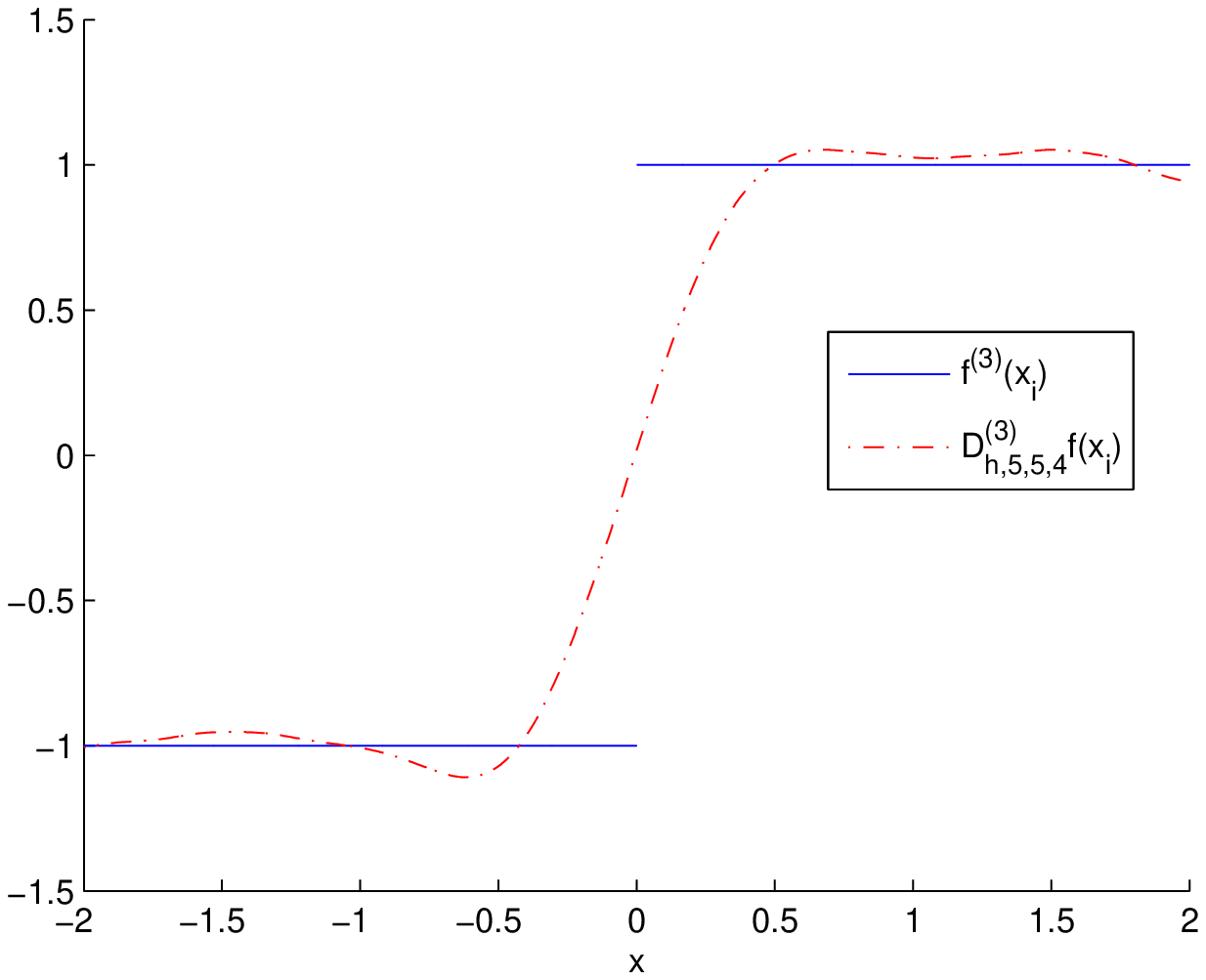}}
  \subfigure[$\delta=0.015$, $n=3$, $\alpha=2$, $q=2$ and $h=1500T_s$.]
 {\includegraphics[scale=0.55]{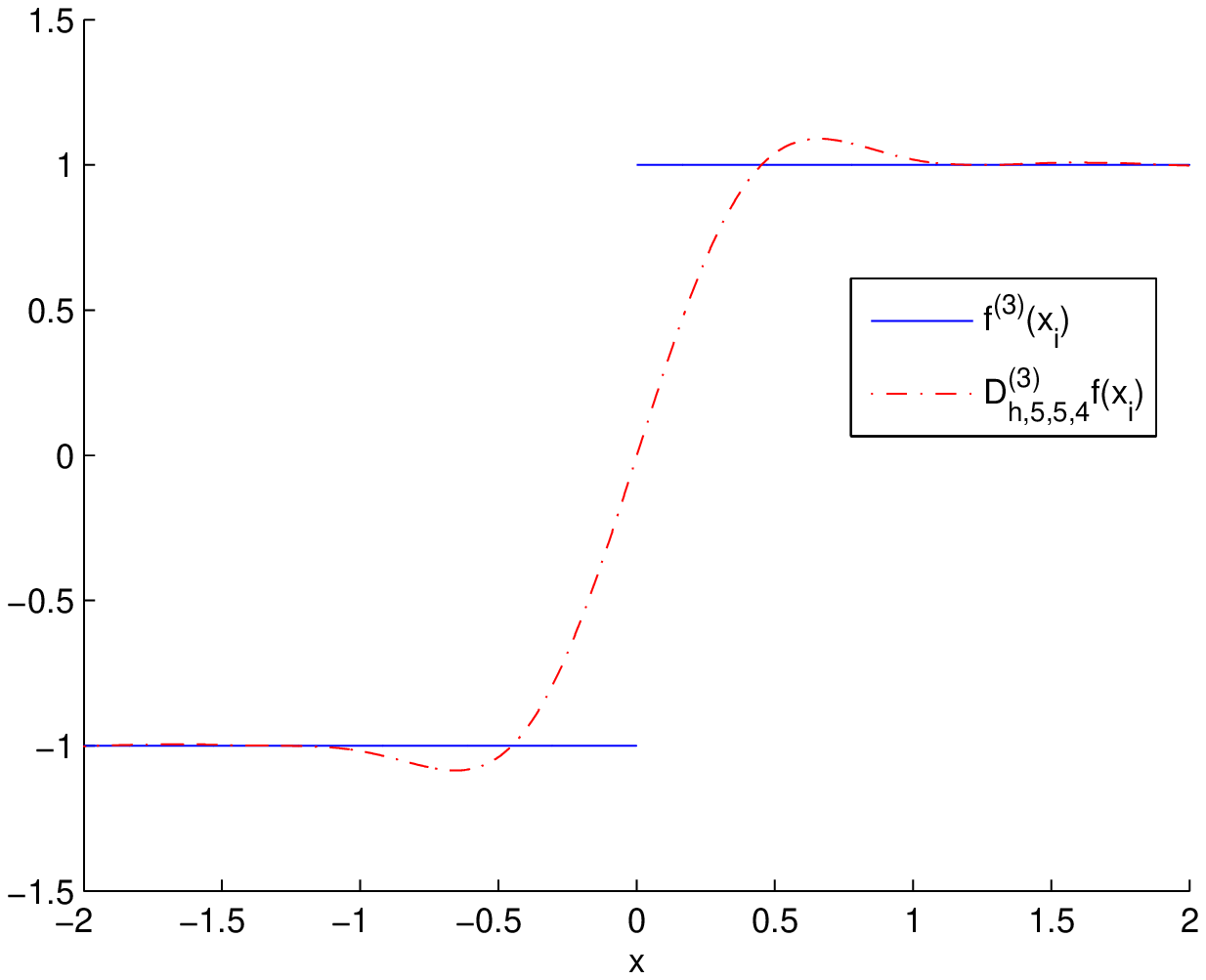}}
\caption{The exact values of  $f_3^{(n)}(x_i)$ and the estimated
values $D_{h,\alpha,\alpha,q}^{(n)}f_3(x_i)$.} \label{fig9}
\end{figure}

\subsection{Simulation results}
The tests are performed by using
Matlab R2007b. Let $f^{\delta}(x_i)= f(x_i) + c \varpi(x_i)$ be
a generated noise data with an equidistant sampling period
$T_s=10^{-3}$ where $c>0$. The
 noise $c \varpi(x_i)$ are simulated from a
zero-mean white Gaussian $iid$ sequence by the Matlab function
'randn' with STATE reset to $0$. By using the well-known three-sigma
rule, we can assume that the noise level for $c \varpi$ is equal to
$3c$. We use the trapezoidal method to approximate the integrals in
our estimators with $2m+1$ values. The estimated derivatives of $f$
at the point $x_i \in I= [-2,2]$  are calculated from the noise data
$f^{\delta}(x_j)$ with $x_j \in [-x_i-h, x_i+h]$, where $h=m T_s$
and $2m+1$ is the number of sampling data used to calculate our
estimation inside the sliding integration windows. When all the
parameters are chosen, $Q_{\alpha,\beta,n,q}$ in the integrals of
our estimators can be calculated explicitly by off-line work with
the $O(n^2)$ complexity. Hence, our estimators can be written like a
discrete convolution product of these pre-calculated coefficients.
Thus, we only need $2m+1$ multiplications and $2m$ additions to
calculate each estimation.

The numerical integration method has an approximation error. Thus,
the total error for our estimators can be bounded by
\begin{equation*}
\begin{split}
\left|T_m \left(Q_{\alpha,n,q}(\cdot)\,f^{\delta}(x_i+h \cdot)\right)- f^{(n)}(x_i) \right| &\leq   \left|T_m \left(Q_{\alpha,n,q}(\cdot)\,f^{\delta}(x_i+h \cdot)\right)- T_m \left(Q_{\alpha,n,q}(\cdot)\,f(x_i+h \cdot)\right) \right| \\ & +
\left|T_m \left(Q_{\alpha,n,q}(\cdot)\,f(x_i+h \cdot)\right)- D_{h,\alpha,\alpha,{q}}^{(n)}f(x_i)\right|+ \left| D_{h,\alpha,\alpha,{q}}^{(n)}f(x_i) - f^{(n)}(x_i)\right|\\
& \leq B_{noise}+ B_{num}+ B_{bias}=B_{total},
\end{split}
\end{equation*}
where $T_m \left(Q_{\alpha,n,q}(\cdot)\,f(x_i+h \cdot)\right)$
(resp. $T_m \left(Q_{\alpha,n,q}(\cdot)\,f^{\delta}(x_i+h
\cdot)\right)$) is the numerical approximation to
$D_{h,\alpha,\alpha,{q}}^{(n)}f(x_i)$ (resp.
$D_{h,\alpha,\alpha,{q}}^{(n)}f^{\delta}(x_i)$) with the trapezoidal
method
  and $B_{num}$ is the well-known error bound for the numerical integration error \cite{Ralston}:
\begin{equation}\label{Nierror}
\left| D_{h,\alpha,\alpha,{q}}^{(n)}f(x_i)- T_m
\left(Q_{\alpha,n,q}(\cdot)\,f(x_i+h \cdot)\right) \right| \leq
\frac{2^3}{12 (2m)^2} \sup_{t\in [-1,1]}
\left(Q_{\alpha,n,q}(t)\,f(x_i+ht)\right)^{(2)} = B_{num}.
\end{equation}

We are going to set the value of $m$ such that $B_{total}$ reaches
its minimum and consequently the total errors in the following two 
examples can be minimized. For this, we need to calculate some values of $f^{(k)}$
with $k=0,\cdots,n+q+2$. According to Remark \ref{remarque3}, we
calculate the value of $M_{n+2+q}$ in the interval
$[-2-\frac{h}{n+q+3},2+\frac{h}{n+q+3}]$. However, in practice, the
function $f$ is
unknown.\\

\noindent\textbf{Example 1.} We choose $f_1(x)=\sin(2\pi x)e^{-x^2}$
as the exact function. The numerical results are shown in Figure
$\ref{fig2}$, where the noise level $\delta$ is equal to $0.15$. The
solid lines represent the exact derivative values of $f_1^{(n)}$ for
$n=1,2,3,4$ and the dash-dotted lines represent the estimated
derivative values $D_{h,\alpha,\alpha,q}^{(n)}f_1(x_i)$. Moreover,
we give in Table $\ref{tab2}$ the total error  values
$\displaystyle\max_{x_i \in [2,2]} \left|
D_{h,\alpha,\alpha,q}^{(n)}f_1(x_i)- f^{(n)}_1(x_i)\right|$ for the
following noise levels: $\delta=0.15$ and  $\delta=0.015$. We can
see also the total error values produced with a larger sampling
period $T'_s=10T_s=10^{-2}$.

\begin{table}[h!]
\caption{$\displaystyle\max_{x_i \in [2,2]} \left|
D_{h,5,5,4}^{(n)}f_1(x_i)- f_1^{(n)}(x_i)\right|$.} \label{tab2}
 \centering

\begin{tabular}{|c|c|c|c|c|}
  \hline
  $\delta$ & $n=1 \ (m)$ & $n=2\ (m)$ & $n=3\ (m)$ & $n=4\ (m)$ \\
  \hline
  $0.15$ & $9.45e-002 \ (591)$ & $1.1\ (698)$ & $1.258e+001\ (777)$ & $1.278e+002 \ (850)$ \\
  \hline
  $0.015$ & $1.85e-002 \ (425)$ & $2.951e-001\ (523)$ & $3.888\ (601)$ & $4.588e+001 \ (675)$ \\
  \hline
 $0.015 \,(T'_s=0.01)$ & $4.06e-002 \ (47)$ & $5.645e-001\ (55)$ & $7.359\ (62)$ & $9.686e+001 \ (69)$ \\
  \hline
\end{tabular}
\end{table}

\noindent\textbf{Example 2.} When $f_2(x)=e^{x^2}$, we give our
numerical results in Figure $\ref{fig1}$ with the noise level
$\delta=0.15$, where the corresponding errors  are given in Figure
$\ref{figerror}$. In  Table $\ref{tab1}$, we also give the total
error values $\displaystyle\max_{x_i \in [2,2]} \left|
D_{h,\alpha,\alpha,q}^{(n)}f_2(x_i)- f_2^{(n)}(x_i)\right|$ for
$\delta=0.15$ and $\delta=0.015$, where  the total error values are
produced with $T_s$ and a larger sampling period $T'_s=10^{-2}$.

\begin{table}[h!]
\caption{$\displaystyle\max_{x_i \in [2,2]} \left|
D_{h,5,5,4}^{(n)}f_2(x_i)- f_2^{(n)}(x_i)\right|$.} \label{tab1}
 \centering

\begin{tabular}{|c|c|c|c|c|}
  \hline
  $\delta$ & $n=1 \ (m)$ & $n=2\ (m)$ & $n=3\ (m)$ & $n=4\ (m)$ \\
  \hline
  $0.15 $ & $1.42e-001\ (442)$ & $2.152\ (549)$ & $2.982e+001\ (643)$ & $3.756e+002\ (733)$ \\
  \hline
  $0.015 $ & $2.22e-002\ (346)$ & $4.435e-001\ (428)$ & $5.973\ (510)$ & $8.769e+001\ (595)$ \\
  \hline
$0.015 \,(T'_s=0.01)$ & $3.404e-001\ (54)$ & $3.425 \ (61)$ & $3.638 e+001\ (68)$ & $5.235e+002\ (79)$ \\
  \hline
\end{tabular}
\end{table}

We can see in Figure $\ref{figerror}$  that  the maximum of the
total error  for each estimation (solid line) is produced nearby the
extremities where the bias term error plus the numerical error
(dash line) are much larger than the noise error.  The noise error
(dash-dotted line) is much larger elsewhere. This is due to the fact
that the total error bound $B_{total}$ is calculated globally in the
interval $[-h-2,2+h]$. The value of $m$ with which $B_{total}$
reaches its minimum is used for all the estimations
$D_{h,\alpha,\alpha,q}^{(n)}f_2(x_i)$ with $x_i \in [-2,2]$. This
value is only appropriate for the estimations nearby the
extremities, but not for the others. In fact, when the bias term
error and the numerical integration error decrease, we should
increase the value of $m$ so as to reduce the noise errors. In order
to improve our estimations, we can choose locally the value of
$m=m_i$, $i.e.$  we search the value $m_i$ which minimizes
$B_{total}$ on $[-h_i+x_i,x_i+h_i]$ where $h_i=m_i T_s$.
 We can see
in Figure \ref{fig11} the errors for these improved estimations
$D_{h_i,5,5,4}^{(n)}f(x_i)$. The different
values of $m_i$ are also given in Figure \ref{fig11}. The
corresponding error bounds are given in Figure \ref{fig22}. We can
observe that the error bounds proposed in this paper are correct but
not optimal. However, the parameters' influence to these error
bounds can help us to know the tendency of errors  so as to choose
parameters for our estimations. On the one hand, the chosen parameters
may not be optimal, but as we have seen in our examples, they give good
estimations. On the other hand, the optimal parameters $q_{op}$,
$\alpha_{op}$ and $m_{op}$ with which the total error bound reaches
its minimum may not give the best estimation. That is why we only use
these error bounds to choose  the value of $m$.

\noindent\textbf{Example 3.} Let us consider the following  function
\begin{equation*}
    f_3(x)=\left\{
       \begin{array}{rl}
         -\displaystyle\frac{1}{6}x^3+2x, & \text{ if } x \leq 0, \\
         \displaystyle\frac{1}{6}x^3+2x, & \text{ if } x > 0,
       \end{array}
     \right.
\end{equation*}
which is $C^{2}$ on $I=[-2,2]$.
The second derivative of $f_3$ is equal to $|x|$. Consequently,  $f_3^{(3)}$ does not
exist at $x=0$. If $n \geq 1$, then this function does not satisfy
the condition $f \in C^{n+2+q}(I)$ of Corollary \ref{corollary5}.
The numerical results are shown in Figure $\ref{fig9}$, where the
sampling period is $T_s=10^{-3}$ and the noise level $\delta$ is equal to
$0.15$ and $0.015$ respectively. The solid lines represent the exact
derivative values of $f_3^{(n)}$ for $n=1,2,3$ and the dash-dotted lines
represent the estimated derivative values
$D_{h,\alpha,\alpha,q}^{(n)}f_3(x_i)$. For the estimations of $f^{(1)}$ and
$f^{(2)}$, we set $\alpha=5$ and $q=4$. When we estimate $f^{(3)}$, the noise error increases. Hence, we need to decrease the values of  $\alpha$ and $q$ to
 $\alpha=2$ and $q=2$. In  Table $\ref{tab3}$, we  give also the total
error values $\displaystyle\max_{x_i \in [2,2]} \left|
D_{h,\alpha,\alpha,q}^{(n)}f_3(x_i)- f_3^{(n)}(x_i)\right|$ for
$n=1,2$ and $\delta=0.015,0.15$.

\begin{table}[h!]
\caption{$\displaystyle\max_{x_i \in [2,2]} \left|
D_{h,5,5,4}^{(n)}f_3(x_i)- f_3^{(n)}(x_i)\right|$.} \label{tab3}
 \centering
\begin{tabular}{|c|c|c|c|c|}
  \hline
  $\delta$ & $n=1 \ (m)$ & $n=2\ (m)$   \\
  \hline
  $0.15 $ & $9.7e-003\ (1700)$ & $9.65 e-002\ (1700)$  \\
  \hline
  $0.015 $ & $4.7e-003\ (1200)$ & $7.23 e-002\ (1200)$ \\
  \hline
\end{tabular}
\end{table}

\end{document}